\documentclass[11pt,amsfonts]{article}
\usepackage{graphicx}
\usepackage{latexsym}
\usepackage{amssymb}
\usepackage{amsmath}
\usepackage{enumerate}

\usepackage{amsmath}
\usepackage{stmaryrd}
\usepackage{hyperref}
\usepackage{amssymb}
\usepackage{CJK}
\usepackage{color}

\usepackage{amsmath,latexsym,amssymb,amsthm,array,amsfonts}

\usepackage{mathrsfs,dsfont,bbm}

\usepackage{layout}
\usepackage{layout}

\usepackage[utf8]{inputenc}
\usepackage[LGR, T1]{fontenc}
\usepackage[english]{babel}
\usepackage{lmodern}
\usepackage{amsmath, esint}
\usepackage{amssymb}
\usepackage{color}
\usepackage{soul}
\usepackage[normalem]{ulem}
\usepackage{enumitem}
\usepackage{graphicx}
\usepackage[usenames, dvipsnames]{xcolor}
\usepackage{stmaryrd}
\usepackage{subfig}
\usepackage{multicol}
\usepackage{multirow}
\usepackage{hyperref}
\usepackage{ifthen}
\usepackage{listings}
\usepackage{diagbox}

\usepackage{tikz}

\newtheorem{prop}{Proposition}
\newtheorem{lemma}{Lemma}

\newtheorem{theorem}{Theorem}

\newtheorem{remark}{Remark}

\def\real{{\mathord{{\rm I\kern-2.8pt R}}}}        
\def\inte{{\mathord{{\rm I\kern-2.8pt N}}}}

\def\sZZ{{\rm Z\kern-2.8ptem{}Z}}

\def\z{{\mathchoice
		{\sZZ}
		{\sZZ}
		{\rm Z\kern-0.30em{}Z}
		{\rm Z\kern-0.25em{}Z} }}
\def\sQQ{{\kern 0.27em \vrule height1.45ex width0.03em depth0em
		\kern-0.30em \rm Q}}
\def\qu{{\mathchoice
		{\sQQ}
		{\sQQ}
		{\kern 0.225em \vrule height1.05ex width0.025em depth0em \kern-0.25em \rm Q}
		{\kern 0.180em \vrule height0.78ex width0.020em depth0em \kern-0.20em \rm Q}
}}
\def\sCC{{\kern 0.27em \vrule height1.45ex width0.03em depth0em
		\kern-0.30em \rm C}}
\def\complex{{\mathchoice
		{\sCC}
		{\sCC}
		{\kern 0.225em \vrule height1.05ex width0.025em depth0em \kern-0.25em \rm C}
		{\kern 0.180em \vrule height0.78ex width0.020em depth0em \kern-0.20em \rm C}
}}

\newcommand{\R}{\mathbb{R}}

\newcommand{\E}{\mathbf{E}}
\renewcommand{\P}{\mathbb{P}}

\def\qed{ \hfill \vrule width.25cm height.25cm depth0cm\smallskip}

\newcommand{\norm}[2]{{\left\Vert #1 \right\Vert}_{#2}}

\newcommand{\ignore}[1]{}
\begin{document}
	
	\renewcommand{\thefootnote}{\fnsymbol{footnote}}
	
	\renewcommand{\thefootnote}{\fnsymbol{footnote}}

	\title{Multidimensional Stein's method for Gamma approximation} 
	
	\author{Ciprian A. Tudor and  J\'er\'emy Zurcher\vspace*{0.2in} \\
	CNRS, Universit\'e de Lille \\
	Laboratoire Paul Painlev\'e UMR 8524\\
	F-59655 Villeneuve d'Ascq, France.\\
	\quad  ciprian.tudor@univ-lille.fr\\
	\quad jeremy.zurcher@univ-lille.fr\\
	\vspace*{0.1in} }
	
	\maketitle

\begin{abstract}
	Let $ F(\nu)$ be the centered Gamma law with parameter $\nu >0$ and let us denote by $\P_{\mathbb{Y}}$ the probability distribution of a random vector $\mathbb{Y}$. We develop a multidimensional variant of the Stein's method for Gamma approximation that allows to obtain bounds for the second Wasserstein distance between the probability distribution of an arbitrary random vector $ (X,  \mathbb{Y})$ in $\mathbb{R}\times \mathbb{R} ^{n}$ and the probability distribution $ F(\nu)\otimes \P_{\mathbb{Y}}$.  In the case of random vectors with components in Wiener chaos, these bounds lead to some interesting criteria for  the joint  convergence of a sequence $ \left( (X_{n}, \mathbb{Y}_{n}), n\geq 1\right)$ to $F(\nu)\otimes \P_{\mathbb{Y}}$, by assuming that $ (X_{n}, n\geq 1)$ converges in law, as $n\to \infty$, to $ F(\nu)$ and $ (\mathbb{Y}_{n}, n\geq 1)$ converges in law, as $n\to \infty$, to an arbitrary random vector $\mathbb{Y}$. We illustrate our criteria by two concrete examples.   
\end{abstract}

	\vskip0.3cm

{\bf 2010 AMS Classification Numbers:}  60F05,60G15,60H05,60H07.

\vskip0.3cm


\section{Introduction}

The Stein's method represents  a popular probabilistic collection of techniques that allows to evaluate the distances between the probability distributions of random variables.  Given two random variables $F, G$, the Stein's method allows to obtain sharp  estimates for the quantities of the form 
\begin{equation}
	\label{intro-1}
	\sup_{h\in \mathcal{H} }\left| \E [h(F)] - \E [h(G)] \right| 
\end{equation}
where $\mathcal{H} $ constitutes a large enough  class of functions. Of particular importance  is the case  when one of the two random variables follows the Gaussian distribution but the cases of other target distributions are also of interest.  We refer to the 
monographs and surveys  \cite{CGS},  \cite{CS}, \cite{Re}, \cite{Stein1} for a detailed description of the techniques of Stein's method and for its applications. 

Our work concerns a variant of the Stein's method that allows to measure the distance between two random vectors with the same marginals but with different correlations between their components. This variant has been initiated in \cite{Pi}, where the author, by combining the Stein's heuristics with the tools of Malliavin calculus, obtained bounds for the Wasserstein distance between the \\ probability distribution of a random vector $(X, Y)$ where $X\sim \mathcal{N}(0,1)$ and $Y$ is an arbitrary random variable and $\P_{Z} \otimes \P_{Y}$, \textit{i.e.} the law of the vector $(Z, Y)$ where $Z \sim \mathcal{N} (0,1)$ and $Z$ is independent by $Y$. We denote by $\P_{\mathbb{X}}$ the probability distribution of a random vector $\mathbb{X}$.  In some sense, this approach allows to evaluate how far are the components $X$ and $Y$ from being independent. The method has been extended in \cite{T3}, by giving asymptotic results and by focusing on the case of random vectors with components in Wiener chaos. 

Our purpose is to develop a similar theory for the case of the centered Gamma distribution with parameter $\nu > 0$, denoted $F (\nu)$ in the sequel. The basic observation is that if $X \sim F(\nu)$ and $X$ is independent from an arbitrary $n$-dimensional random vector $\mathbb{Y}$, then 
\begin{equation}
	\label{intro-2}
	\E \left[  2(X+\nu) \frac{\partial f}{\partial x} (X, \mathbb{Y})- X f(X, \mathbb{Y})\right] = 0
\end{equation}
for a large class of differentiable functions $ f: \mathbb{R}\times  \mathbb{R} ^{n} \to \mathbb{R}$.  Also, if (\ref{intro-2}) holds  true for a large class of functions $f$, then $X\sim F(\nu)$ and $X$ is independent of $\mathbb{Y}$. Then, we introduce a multidimensional counterpart of the standard Sten's equation for the Gamma law, i.e.
\begin{equation}
	\label{intro-3}
	2(x+\nu) \frac{\partial f}{\partial x}(x, {\bf y}) -xf(x, {\bf y}) = h(x,\mathbf{y}) - \E [ h(Z_{\nu}, \mathbf{y}) ],
	\end{equation}
where $Z_{\nu } \sim F(\nu)$ and $h: \mathbb{R}\times \mathbb{R} ^{n} \to \mathbb{R}$ belongs to suitable class of functions.  We analyze in details the existence, the uniqueness and the regularity of the solution to (\ref{intro-3}) and of its partial derivatives. In particular, we prove that  there is a unique bounded solution to (\ref{intro-3}) and  the infinity norm of this solution and of its first order  partial derivatives are controlled by the infinity norm of $h$ and of its first and second order partial derivatives. By combining the Stein's equation with the Malliavin's integration by parts, we obtain bounds for 
\begin{equation*}
	d_{2} \left( \P_{(X,\mathbb{Y})}, \P_{Z_{\nu}} \otimes \P_{\mathbb{Y}}\right),
\end{equation*}
where$d_{2}$ is the second Wasserstein distance (defined later in our work),  $Z_{\nu}\sim F(\nu)$, and $X, \mathbb{Y}$ are arbitrary random vectors in $\mathbb{R}$ and $\mathbb{R}^{n}$ respectively, sufficiently regular in the sense of Malliavin calculus. We also  analyze the case of random sequences in Wiener chaos where those general results translate into easy-to-apply criteria for Gamma approximation. We consider a sequence $ (X_{k}, k\geqslant 1)$ and a sequence of $n$-dimensional random vectors $ (\mathbb{Y}_{k}, k\geqslant 1)$ such that for each $k\geqslant 1$, $X_{k}$ and the components of $\mathbb{Y}_{k}$ belongs to a Wiener chaos and we assume that, in distribution, $X_{k}$ converges to $Z_{\nu} \sim F(\nu)$ and $\mathbb{Y}_{k}$ converges to an arbitrary random vector $\mathbb{Y}$ as $k\to \infty$. Under some rather natural conditions (and easy to check in particular case), we prove that the vector $ \left( (X_{k}, \mathbb{Y}_{k}), k\geqslant 1\right)$ converges in law, as $k\to \infty$, to $F(\nu) \otimes \P_{\mathbb{Y}}$. We also evaluate the corresponding rate of convergence  under the so-called second Wasserstein distance.    The method is illustrated by some examples. Our results extend the methods and findings in the literature related to the Gamma approximation, see e.g. \cite{AAPS}, \cite{AAPS2}, \cite{AS}, \cite{AET}, \cite{APP}, \cite{DP}, \cite{D}, \cite{KuTu}, \cite{KuTu2}, \cite{NP1}, \cite{NP2}.

We organized our paper as follows.  In Section 2 we give a characterization for the probability distribution of a random vector whose first marginal is a centered Gamma law and the rest of the vector is independent by the first component. This naturally leads to a multidimensional Stein's equation corresponding to the law of a such multidimensional random vectors. We then analyze in details the solution to the Stein's equation on the whole real line and we prove suitable bounds for this solution and for its partial derivatives. In Section 3 we combine our findings on the solution to the Stein's equation with the techniques of Malliavin calculus, in the spirit of \cite{NP1}, in order to obtain bounds for the second Wasserstein distance between an arbitrary random vector and a random vector with the first component following the $F(\nu)$ law and the rest of the vector independent by this first component. Section 4 is devoted to the study of sequences of random variables in Wiener chaos. If $(X_{n}, n\geq 1)$ converges in distribution to $ F(\nu)$ and $ (\mathbb{Y}_{n}, n\geq 1)$ converges in distribution to an arbitrary law $U$  on $\mathbb{R}^{n}$, we give criteria to check the  joint convergence of the vector $((X_{n}, \mathbb{Y}_{n}), n\geq 1)$ to $ F(\nu)\otimes U$. We illustrate our result by two concrete examples in Section 5. Finally, the Appendix contains the basics of Malliavin calculus and the proof of as technical result. 

\section{The multidimensional Stein's equation and its solution}
Let $ Z_{\nu}$ be a random variable following the centered Gamma distribution $F(\nu)$ with parameter $\nu >0$, whose probability density function is 
\begin{equation}\label{pnu}
	p_{\nu} (x) := \frac{ (x+\nu) ^{\frac{\nu}{2}-1} e ^{-\frac{x+\nu}{2}}}{2 ^{\frac{\nu}{2}}\Gamma \left( \frac{\nu}{2}\right)} \mathbf{1}_{ (-\nu, + \infty)} (x).
\end{equation}
When $\nu $ is an integer, then $Z_{\nu}$ coincides in law with $\sum_{i=1}^{\nu} (N_{i}^{2}-1)$, where $N_{1},...,N_{\nu}$ are standard Gaussian independent random variables. We recall that in the standard Stein's method for Gamma approximation, the development of the method is based by the fact that $ Z_{\nu}\sim F(\nu)$ if and only if 
\begin{equation}\label{s1}
	\E \left[  2 (Z_{\nu}+\nu)_{+} f'(Z_{\nu})\right]  = \E \left[ Z_{\nu} f(Z_{\nu})\right]
\end{equation}
for a large class of differentiable functions $f: \mathbb{R} \to \mathbb{R}$, where $x_{+}=\max(x,0)$. Let us state and prove an analogous result in the multidimensional context.  We recall that if $X\sim F(\nu)$, then its characteristic function $ \psi $ is the unique solution to the differential equation
\begin{equation}\label{pde}
	\forall \lambda \in \R, \hskip 0.3cm (1-2i\lambda) \psi'(\lambda)+ 2 \lambda \nu \psi (\lambda)=0,  
\end{equation}
with $\psi (0)=1$.

\begin{lemma} Let $\mathbb{Y}=(Y_{1},...,Y_{n})$ be a $n$-dimensional real random vector. 
	\begin{enumerate} 
		\item Assume $ Z_{\nu }\sim  F(\nu ) $ and $Z_{\nu}$ and $\mathbb{Y}$ are independent. Then 
		\begin{equation*}
			\E\left[  2 ( Z_{\nu }+\nu) \frac{\partial f}{\partial x} (Z_{\nu}, \mathbb{Y})\right] = \E\left[ Z_{\nu} f( Z_{\nu}, \mathbb{Y})\right] ,
		\end{equation*}
	\end{enumerate}
	\noindent for any function $ f: \mathbb{R}\times \mathbb{R} ^{n} \to \mathbb{R} $ such that $\E \left[ \left|  (Z_{\nu }+\nu) \frac{\partial f}{\partial x} (Z_{\nu}, \mathbb{Y})\right| \right] <\infty$ and $ \E \left[ \left| Z_{\nu} f( Z_{\nu}, \mathbb{Y})\right| \right]  <\infty$. 
	\begin{enumerate}[resume]
		\item Conversely, assume that $ Z$ is an integrable random variable such that
			\begin{equation}\label{23i-1}
			\E \left[ 2 ( Z+\nu)_{+} \frac{\partial f}{\partial x} (Z, \mathbb{Y})\right] = \E \left[ Zf( Z, \mathbb{Y})\right],
		\end{equation}
    \end{enumerate}
		\noindent for any function $ f: \mathbb{R}\times \mathbb{R}^{n} \to \mathbb{R} $ such that $\E \left[ \left|  (Z+\nu) \frac{\partial f}{\partial x} (Z, \mathbb{Y})\right| \right] < \infty$ and $ \E \left[ \left| Zf( Z, \mathbb{Y})\right| \right]  <\infty$. Then $ Z \sim F (\nu)$ and $ Z $ is independent of $\mathbb{Y}$. 
	
\end{lemma}
\noindent {\bf Proof: } Assume that $ Z_{\nu} \sim F(\nu) $ and $Z_{\nu}$ and $\mathbb{Y}$ are independent. Then, by the Stein's characterization of the centered Gamma law  (\ref{s1}), we have for every ${\bf y}\in \mathbb{R }^{n}$, 

\begin{equation*}
	\E \left[ 2 (Z_{\nu}+\nu) \frac{\partial f}{\partial x} (Z_{\nu}, {\bf y})\right] = \E \left[  Z_{\nu} f (Z_{\nu}, {\bf y})\right].
\end{equation*}
So, 
\begin{equation*}
	\int_{-\nu}^{\infty} 2 (x+ \nu) \frac{\partial f}{\partial x} (x, {\bf y}) \textrm{d} \P_{Z_{\nu}}(x)= \int_{-\nu}^{\infty} x f(x, {\bf y}) \textrm{d} \P_{Z_{\nu}} (x),
\end{equation*}
for every ${\bf y}\in \mathbb{R} ^{n}$.  We integrate the above relation with respect to $\P_{\mathbb{Y}}$ on $\mathbb{R}^{n}$ and we find for the left-hand side, by using the independence of $ Z_{\nu}$ and $\mathbb{Y}$,
\begin{eqnarray*}
	\int_{\mathbb{R} ^{n}} \left( 	\int_{-\nu}^{\infty} 2 (x+ \nu) \frac{\partial f}{\partial x} (x, {\bf y}) \textrm{d} \P_{Z_{\nu}}(x)\right) \textrm{d} \P_{\mathbb{Y}}({\bf y})
		&=&\int_{ (-\nu, \infty)\times \mathbb{R} ^{n}}  2 (x+ \nu) \frac{\partial f}{\partial x} (x, {\bf y}) \textrm{d}(\P_{Z_{\nu}} \otimes \P_{\mathbb{Y}}) (x, {\bf y})\\
		&=& \E \left[   2 (Z_{\nu}+\nu) \frac{\partial f}{\partial x} (Z_{\nu}, \mathbb{Y}) \right] .
\end{eqnarray*}
Similarly for the right-hand side 
\begin{eqnarray*}
	\int_{\mathbb{R} ^{n}} \left( \int_{-\nu}^{\infty} x f(x, {\bf y}) \textrm{d} \P_{Z_{\nu}}(x) \right) \textrm{d} \P_{\mathbb{Y}}({\bf y})&=& \int_{ (-\nu, \infty)\times \mathbb{R} ^{n}}  x f(x, {\bf y}) \textrm{d}(\P_{Z_{\nu}} \otimes \P_{\mathbb{Y}}) (x, {\bf y})\\
	&=& \E \left[ Z_{\nu } f(Z_{\nu}, \mathbb{Y}) \right] . 
\end{eqnarray*}
Let us now prove the second point of the lemma. By taking a function $f$ with support contained in $(-\infty, -\nu) $, we found that $ \P (Z_{\nu}\leq -\nu )=0$. Let $\varphi$ be the characteristic function of the random vector $(Z, \mathbb{Y})$, \textit{i.e.}
\begin{equation*}
	\varphi( \lambda_{1},  \lambda _{2}) := \E \left[ e ^{i(\lambda _{1}Z+ \lambda _{2} \cdot \mathbb{Y})}\right],
\end{equation*}
for $\lambda _{1} \in \mathbb{R}, \lambda _{2} \in \mathbb{R} ^{n}$, where ${\bf x}\cdot {\bf y}$ denotes the Euclidean scalar product  of ${\bf x}, {\bf y}\in \mathbb{R} ^{n}$. We take the derivative of $\varphi$ with respect to $\lambda _{1}$ and we apply (\ref{23i-1}) (for the real and imaginary parts of $\varphi$). By denoting $g(x, {\bf y})= e ^{i(\lambda _{1}x+\lambda _{2}\cdot  {\bf y})}$, then since $Z$ is integrable, $g$ satisfies the conditions to have (\ref{23i-1}) and we get 
\begin{eqnarray*}
\frac{\partial \varphi}{\partial \lambda _{1}} (\lambda _{1}, \lambda _{2})  & = &	 i \E \left[ Z e ^{i(\lambda _{1} Z+ \lambda _{2} \cdot \mathbb{Y})} \right] = i \E [ Z g(X, \mathbb{Y})] \\
&=& i \E \left[ 2 (Z+\nu) \frac{\partial g}{\partial x} (Z, \mathbb{Y})\right] = -2 \lambda _{1} \E \left[ (Z+\nu)e ^{i(\lambda _{1} Z+ \lambda _{2} \cdot \mathbb{Y})} \right] \\
&=& 2i \lambda _{1} \frac{\partial \varphi}{\partial \lambda _{1}}  (\lambda _{1}, \lambda _{2}) - 2\nu \lambda _{1} \varphi (\lambda _{1}, \lambda _{2}).
\end{eqnarray*}
Consequently, for every $\lambda _{1} \in \mathbb{R}, \lambda _{2} \in \mathbb{R} ^{n}$, we have 
\begin{equation*}
	(1-2i\lambda _{1} )\frac{\partial \varphi}{\partial \lambda _{1}}  (\lambda _{1}, \lambda _{2}) +2\lambda _{1} \nu \varphi (\lambda_{1}, \lambda _{2})=0.
\end{equation*}
By noticing that for every $\lambda_{2} \in \mathbb{R} ^{n}$
\begin{eqnarray*}
	\varphi(0, \lambda _{2})= \E \left[ e ^{i \lambda _{2} \cdot \mathbb{Y}} \right] = \varphi_{\mathbb{Y}}(\lambda _{2}),
\end{eqnarray*}
where $\varphi_{\mathbb{Y}}$ stands for the characteristic function of the random vector $\mathbb{Y}$, we obtain from (\ref{pde}),  for every $\lambda _{1} \in \mathbb{R}, \lambda _{2} \in \mathbb{R} ^{n}$,
\begin{equation*}
	\varphi(\lambda _{1}, \lambda _{2})= \psi (\lambda _{1}) \varphi_{\mathbb{Y}}(\lambda _{2}).
\end{equation*}
This means that $ Z\sim F(\nu)$ and $Z$ is independent of $\mathbb{Y}$. \qed

Let us introduce the multidimensional Stein's equation for the centered Gamma law $F(\nu)$
\begin{equation}
	\label{s2}2(x+\nu) \frac{\partial f}{\partial x} (x, {\bf y})- xf(x, {\bf y})= h(x,y) - \E \left[ h(Z_{\nu}, {\bf y}) \right],
\end{equation}
with $Z_{\nu}\sim F(\nu) $ and $h: \mathbb{R} \times \mathbb{R} ^{n} \to \mathbb{R}$ such that $\E \left[ | h(F(\nu), {\bf y}) | \right] <\infty$ for any ${\bf y}\in \mathbb{R} ^{n}$. In order to solve (\ref{s2}), let us recall some facts concerning the one-dimensional Stein's equation (on the whole real line) for the Gamma law
\begin{equation}
	\label{se1}
	2(x+\nu)f'(x)-xf(x)= h(x) - \E[h(F(\nu))],
\end{equation} 
with $h:\mathbb{R}\to \mathbb{R}$ satisfying $\E \left[ \left| h(F(\nu)) \right| \right] <\infty$. Recall that the density of this law is given by (\ref{pnu}) and let us also introduce the function
\begin{equation}{\label{q}}
	q_{\nu} (x) := \frac{(-(x+ \nu))^{\frac{\nu}{2} - 1} e^{- \frac{x + \nu}{2}}}{2^{\frac{\nu}{2}} \Gamma \left( \frac{\nu}{2} \right)} \mathbf{1}_{(- \infty, - \nu)} (x).
\end{equation}

Both functions $p_{\nu}$ and $q_{\nu}$ are positive. Note that, contrary to $p_{\nu}$, the function $q_{\nu}$ does not define a probability measure. We consider $F_{\nu}$ the cumulative probability function of $p_{\nu}$. We denote by $\tilde{F}_{\nu}$ the following function defined on $(- \infty, - \nu]$:

\begin{equation} \label{23i-2}
	 \tilde{F}_{\nu} (x) = \int_{x}^{- \nu} q_{\nu} (u) \ \textrm{d}u. 
	\end{equation}
Then $\tilde{F}_{\nu}$ is a positive decreasing map on $(- \infty, - \nu]$.

\begin{lemma}\label{ll2}
	Let $h: \mathbb{R} \to \mathbb{R}$ be measurable such that $\E \left[ \vert h(Z_{\nu})\vert \right]  <\infty$. Then the Stein's equation (\ref{se1}) admits a unique bounded solution given by 
	\begin{eqnarray}
		f_{h}(x) = & & \left[ \int_{-\nu} ^{x} \frac{p_{\nu}(u)}{2(x+\nu)p_{\nu}(x)}\left( h(u)- \E \left[ h(Z_{\nu}) \right] \right) \ \mathrm{d}u \right] \mathbf{1}_{(-\nu, + \infty)}(x) \nonumber \\
		& - & \left[ \int_{x} ^{-\nu}  \frac{q_{\nu}(u)}{2(x+\nu) q_{\nu}(x)}\left( h(u) - \E [ h(Z_{\nu}) ] \right) \ \mathrm{d}u \right] \mathbf{1}_{(-\infty, -\nu)}(x) \label{sol1} \\
		& + & \frac{h(-\nu)- \E \left[ h(Z_{\nu}) \right]}{\nu} \mathbf{1}_{\{ - \nu \}} (x) \nonumber.
	\end{eqnarray}
\end{lemma}

\noindent {\bf Proof: } See Section 2.2 in \cite{DP}. \qed

\begin{remark}\label{rem1}
We can show that this unique solution is $C^1$ on $\R$, by using the l'Hôpital's rule several times, when $h$ is $C^1$ with $h'$ absolutely continuous. Indeed, for instance, for $x$ going to $- \nu$ from above :

\begin{eqnarray*}
    \lim_{x \downarrow -\nu} f'_h (x) = & & \lim_{x \downarrow - \nu} \frac{2(x+ \nu) p_{\nu}(x) \tilde{h}(x) + x \int_{-\nu}^x \tilde{h}(u) p_{\nu}(u) \ \mathrm{d}u}{4(x+\nu)^2 p_\nu (x)} \\
    \overset{(\hat{H})}{=} & & \lim_{x \downarrow - \nu} \frac{2(x+\nu) p_{\nu} (x) h'(x) - x p_{\nu} (x) \tilde{h}(x) + \int_{- \nu}^x \tilde{h}(u) p_{\nu} (u) \ \mathrm{d}u + x \tilde{h}(x) p_{\nu}(x)}{2(2-x)(x+\nu) p_{\nu}(x)} \\
    = & & \frac{1}{2+\nu} \lim_{x \downarrow - \nu} \frac{2(x+\nu)p_{\nu} (x) h'(x) + \int_{- \nu}^x \tilde{h}(u) p_{\nu} (u) \ \mathrm{d}u}{2(x+\nu) p_{\nu} (x)} \\
    \overset{(\hat{H})}{=} & & \frac{1}{2+\nu} \lim_{x \downarrow - \nu} \frac{2 h''(x) (x+ \nu) p_{\nu}(x) - x p_{\nu} (x) h'(x) + \tilde{h}(x) p_{\nu}(x)}{-x p_{\nu}(x)} \\
    = & & \frac{h'(-\nu)}{2+\nu} + \frac{\tilde{h}(-\nu)}{\nu (2 + \nu)},
\end{eqnarray*}
where $\tilde{h}(x) = h(x) - \E[h(Z_{\nu})]$, "$\overset{(\hat{H})}{=}$" means that we used the l'Hôpital's rule and we also applied  the fundamental property 
\begin{equation}
	\label{5mm-1}
2(x+\nu) p_{\nu} (x) = - \int_{- \nu}^x u p_{\nu} (u) \ \mathrm{d}u.
\end{equation}
This works analogously  when $x$ goes to $- \nu$ from below, via a similar relation for $q_{\nu}$ : 
\begin{equation*}
2(x+\nu) q_{\nu} (x) = \int_{x}^{- \nu} u q_{\nu} (u) \ \mathrm{d}u  \mbox{ for every }x < - \nu.
\end{equation*} 

\end{remark}

In a first step, we deduce the existence of a unique bounded solution of the multidimensional Stein's equation (\ref{s2}). 

\begin{prop}
	Let $h: \mathbb{R} \times \mathbb{R} ^{n} \to \mathbb{R}$ be a test function  such that $\E \left[ | h(F(\nu), {\bf y}) | \right] <\infty$ for any ${\bf y}\in \mathbb{R} ^{n}$.  Then (\ref{s2}) admits an unique bounded solution which can be expressed as 
	\begin{eqnarray}
		f_h (x, \mathbf{y}) = & & \left[ \int_{- \nu}^{x} \frac{p_{\nu} (u)}{2(x+\nu) p_{\nu} (x)}  \left( h(u, \mathbf{y}) - \E \left[ h(Z_{\nu}, \mathbf{y}) \right] \right) \ \mathrm{d}u \right] \mathbf{1}_{(-\nu, + \infty)} (x) \nonumber \\
		& - & \left[ \int_{x}^{- \nu} \frac{q_{\nu} (u)}{2(x+\nu) q_{\nu} (x)}  \left( h(u, \mathbf{y}) - \E \left[ h(Z_{\nu}, \mathbf{y}) \right] \right) \ \mathrm{d}u \right] \mathbf{1}_{(-\infty, - \nu)} (x),\label{fh} \\
		& + & \frac{h(-\nu, \mathbf{y})- \E \left[ h(Z_{\nu}, \mathbf{y}) \right]}{\nu} \mathbf{1}_{\{ - \nu \}} (x) \nonumber.
	\end{eqnarray}
\end{prop}

\noindent{\bf Proof: } Let us prove first that (\ref{fh}) satisfies (\ref{s2}). For $x > - \nu$, we have by differentiating (\ref{fh}) with respect to $x$:

\begin{equation*}
    \frac{\partial f_h}{\partial x} (x, \mathbf{y}) = \frac{h(x, \mathbf{y}) - \E[h(Z_{\nu}, \mathbf{y})]}{2 (x + \nu)} + \frac{x}{4(x+\nu)^2 p_{\nu} (x) } \int_{- \nu}^x p_{\nu} (u) (h(x, \mathbf{y}) - \E[h(Z_{\nu}, \mathbf{y})]) \ \textrm{d}u,
\end{equation*}

\noindent and so we conclude that we have (\ref{s2}) in $(-\nu, + \infty)$. For $x < - \nu$ :

\begin{equation*}
    \frac{\partial f_h}{\partial x} (x, \mathbf{y}) = \frac{h(x, \mathbf{y}) - \E[h(Z_{\nu}, \mathbf{y})]}{2 (x + \nu)} - \frac{x}{4(x+\nu)^2 q_{\nu} (x) } \int_x^{- \nu} q_{\nu} (u) (h(x, \mathbf{y}) - \E[h(Z_{\nu}, \mathbf{y})]) \ \textrm{d}u,
\end{equation*}

\noindent giving the same conclusion on $(- \infty, - \nu)$. Finally, we can check by l'Hôpital's rule, in the same way as in \cite{D}, that when $x$ goes to $- \nu$, the solution (\ref{fh}) goes to $\frac{h(-\nu,\mathbf{y})-\E \left[ h(Z_{\nu}, \mathbf{y}) \right]}{\nu}$, either the taking the lft or the right limit. Hence, $x \longmapsto f_h(x, \mathbf{y})$ is continuous on $\R$. Moreover, we can prove, when $h$ is $C^1$ with $\frac{\partial h}{\partial x}$ absolutely continuous, by the same argument as in Remark \ref{rem1}, that is also the case of $\frac{\partial f_h}{\partial x}$. So $f$ is $C^1$ with respect to $x$. This concludes the existence part.

Let us prove the unicity. Assume $f_{h}, g_{h}$ are two bounded solutions to (\ref{s2}). Then for every $(x, {\bf y})\in \mathbb{R} \times \mathbb{R} ^{n}$,
\begin{equation*}
	2(x+\nu) \frac{ \partial (f_{h}-g_{h})}{\partial x}  (x, {\bf y}) -x (f_{h}-g_{h})(x, {\bf y}) =0.
\end{equation*}
By solving the above equations, we get
\begin{equation*}
	(f_{h}- g_{h})(x,{\bf y})= \begin{cases}
		c_{1}({\bf y})\frac{ e ^{\frac{x}{2}}}{(x+\nu) ^{\frac{\nu}{2}}} \mbox{ if } (x, {\bf y})\in (-\nu, \infty)\times \mathbb{R} ^{n}\\
		c_{2}({\bf y}) \frac{ e ^{\frac{x}{2}}}{(-(x+\nu)) ^{\frac{\nu}{2}}} \mbox{ if } (x,{\bf y})\in (-\infty, -\nu)\times \mathbb{R} ^{n}, 
	\end{cases}
\end{equation*}
so
\begin{equation*}
	f_{h} (x, {\bf y}) =  \begin{cases}
	g_{h}(x, {\bf y})+	c_{1}({\bf y})\frac{ e ^{\frac{x}{2}}}{(x+\nu) ^{\frac{\nu}{2}}} \mbox{ if } (x, {\bf y})\in (-\nu, \infty)\times \mathbb{R} ^{n}\\
	g_{h}(x, {\bf y})+	c_{2}({\bf y}) \frac{ e ^{\frac{x}{2}}}{(-(x+\nu)) ^{\frac{\nu}{2}}} \mbox{ if } (x,{\bf y})\in (-\infty, -\nu)\times \mathbb{R} ^{n}.
\end{cases}
\end{equation*}
Consequently,  $f_{h}$ is bounded if and only if $ c_{1}({\bf y})= c_{2}({\bf y}) =0 $ for every $\mathbf{y} \in \mathbb{R}^{n}$. \qed 

 The next step is to find suitable estimates  for the solution (\ref{fh}) and for its partial derivatives.  We follow the methodology proposed in \cite{DP} and \cite{D} in the one-dimensional case. Let us start by an useful result that gives an alternative expression for the solution to (\ref{s2}).  The proof of the below lemma can be found in the Appendix.
 
 \begin{lemma}\label{ll3}
 	
 	\begin{enumerate}
 		\item Let $f_{h}$ be given by (\ref{fh}). Then for $x>-\nu$ and ${\bf y}\in \mathbb{R} ^{n}$, 
 		\begin{eqnarray}
 			& & f_h (x, \mathbf{y}) \\ {\label{fhfub+}}   
 			 = & &  \frac{1 - F_{\nu}(x)}{ \int_{-\nu}^x u p_{\nu} (u) \ \mathrm{d}u} \int_{- \nu}^{x} \frac{\partial h}{\partial x} (w, \mathbf{y}) F_{\nu} (w) \ \mathrm{d}w + \frac{F_{\nu}(x)}{ \int_{- \nu}^{x} u p_{\nu} (u) \ \mathrm{d}u} \int_{x}^{+ \infty} \frac{\partial h}{\partial x} (w, \mathbf{y}) (1 - F_{\nu} (w)) \ \mathrm{d}w,  \nonumber
 		\end{eqnarray}
\noindent and for  every $x < - \nu$ and ${\bf y}\in \mathbb{R} ^{n}$, we have:
 		\begin{eqnarray}
 		f_h (x, \mathbf{y})	= & \displaystyle \frac{1}{\int_x^{- \nu} u q_{\nu} (u) \ \mathrm{d}u} & \left[ \int_{x}^{- \nu} \frac{\partial h}{\partial x} (w, \mathbf{y}) \left[\tilde{F}_{\nu} (x) - \tilde{F}_{\nu} (w) \right] \ \mathrm{d}w \right. \nonumber \\ 
 			& & \left. + \tilde{F}_{\nu}(x)  \int_{-\nu}^{+ \infty} \frac{\partial h}{\partial x} (w, \mathbf{y}) (1 - F_{\nu} (w)) \ \mathrm{d}w \right]. {\label{fhfub-}}
 		\end{eqnarray}

 		\item If  $F_{\nu}$ denotes the cumulative distribution function of the centered Gamma law $F(\nu)$, then for any  $x \in (-\nu, \infty)$ :
 		
 		$$ \int_{- \nu}^x F_{\nu} (s) \ \mathrm{d}s = x F_{\nu} (x) - \int_{- \nu}^x s p_{\nu}(s) \ \mathrm{d}s. $$
 		and 
 		
 		$$ \int_x^{+ \infty} (1 - F_{\nu} (s)) \ \mathrm{d}s = -x(1 - F_{\nu} (x)) + \int_{x}^{+ \infty} s p_{\nu}(s) \ \mathrm{d}s. $$
 		
 		\item  Let $\tilde{F}_{\nu}$ be given by (\ref{23i-2}). Then we have for $x\in (-\infty, -\nu)$,
 		
 		$$ \int_x^{- \nu} \tilde{F}_{\nu} (w) \ \mathrm{d}w = \int_x^{- \nu} u q_{\nu} (u) \ \mathrm{d}u - x \tilde{F}_{\nu} (x). $$
 		
 	\end{enumerate}
 	
 \end{lemma}

Let us now obtain some estimates for the solution to the Stein's equation.  The idea of the proof is based on the proof of Theorem 3.13 in \cite{D}.

\begin{prop}\label{pp2}
	Let $h: \mathbb{R} \times \mathbb{R} ^{n}$ be a $ C^{2}$ map with bounded partial derivatives. Let $ f_{h}$ be the unique bounded solution to (\ref{s2}) corresponding to $h$. We have the following estimates :
	
	\begin{equation}{\label{ingfh}}
		\norm{f_h}{\infty} \leqslant C\norm{\frac{\partial h}{\partial x}}{\infty},
	\end{equation}
	\noindent and for every $j \in \{ 1,...,n \}$, 
	\begin{equation}{\label{ingdfh}}
		\norm{\frac{\partial f_h}{\partial y_j}}{\infty} \leqslant  C\norm{\frac{\partial^2 h}{\partial x \partial y_j}}{\infty} \quad  \text{ and } \quad \norm{\frac{\partial f_h}{\partial x}}{\infty} \leqslant C \norm{\frac{\partial h}{\partial x}}{\infty} .
	\end{equation}
	\end{prop}
\noindent {\bf Proof: } Let us prove the estimation (\ref{ingfh}). Let $x > - \nu$. By using Lemma \ref{ll3}, we have 

\begin{eqnarray*}
	& & |f_h (x, \mathbf{y})| \\ 
	= & & \left| \frac{F_{\nu} (x)}{\int_{- \nu}^x u p_{\nu} (u) \ \textrm{d}u} \int_x^{+ \infty} \frac{\partial h}{\partial x} (w, \mathbf{y}) (1 - F_{\nu} (w)) \ \textrm{d}w + \frac{1 - F_{\nu} (x)}{\int_{- \nu}^x u p_{\nu} (u) \ \textrm{d}u} \int_{- \nu}^x \frac{\partial h}{\partial x} (w, \mathbf{y}) F_{\nu} (w) \ \textrm{d}w \right|.
\end{eqnarray*}

\noindent Note that $\int_{- \nu}^x u p_{\nu} (u) \ \textrm{d}u = - 2 (x+\nu) p_{\nu} (x) \leqslant 0$. Then,  via point 2. in Lemma \ref{ll3}, 

\begin{eqnarray*}
	|f_h (x, \mathbf{y})| \leqslant & \displaystyle \frac{\norm{\frac{\partial h}{\partial x}}{\infty}}{- \int_{- \nu}^x u p_{\nu} (u) \ \textrm{d}u} &  \left( F_{\nu} (x) \int_x^{+ \infty} (1 - F_{\nu}(w)) \ \textrm{d}w +  (1 - F_{\nu} (x)) \int_{- \nu}^x F_{\nu}(w) \ \textrm{d}w \right) \\
= & \displaystyle \frac{\norm{\frac{\partial h}{\partial x}}{\infty}}{- \int_{- \nu}^x u p_{\nu} (u) \ \textrm{d}u} &  \left( F_{\nu} (x) \left( \int_x^{+ \infty} s p_{\nu} (s) \ \textrm{d}s - x(1 - F_{\nu} (x)) \right) \right. \\ 
& & \left. +  (1 - F_{\nu} (x)) \left( x F_{\nu} (x) - \int_{- \nu}^x s p_{\nu} (s) \ \textrm{d}s  \right) \right) \\
	= & \displaystyle \frac{\norm{\frac{\partial h}{\partial x}}{\infty}}{- \int_{- \nu}^x u p_{\nu} (u) \ \textrm{d}u} &  \left( F_{\nu} (x) \int_{-\nu}^{+ \infty} s p_{\nu} (s) \ \textrm{d}s - \int_{- \nu}^x s p_{\nu} (s) \ \textrm{d}s \right).
\end{eqnarray*}
\noindent By using the fact that $\int_{- \nu}^{+ \infty} s p_{\nu} (s) \ \textrm{d}s = \E[Z_{\nu}] = 0$, we conclude the estimation (\ref{ingfh}) on $(- \nu, + \infty)$. 

Now, we proceed in the same way for $x < - \nu$. By Lemma \ref{ll3}

\begin{eqnarray*}
	|f_h (x, {\bf y})| = & \displaystyle \frac{1}{- \int_x^{- \nu} u q_{\nu} (u) \ \mathrm{d}u} & \left|  \left[ \int_{x}^{- \nu} \frac{\partial h}{\partial x} (w, \mathbf{y}) \left[\tilde{F}_{\nu} (x) - \tilde{F}_{\nu} (w) \right] \ \mathrm{d}w \right. \right. \\ 
	&  & \left. \left. + \tilde{F}_{\nu}(x) \int_{- \nu}^{+ \infty} \frac{\partial h}{\partial x} (w, \mathbf{y}) (1 - F_{\nu} (w)) \ \mathrm{d}w \right] \right|.
\end{eqnarray*}

\noindent Since $\tilde{F}_{\nu}$ is decreasing on $(- \infty, - \nu)$ and since $\int_{x}^{- \nu} u q_{\nu} (u) \ \textrm{d}u = 2(x+\nu) q_{\nu} (x) \leqslant 0$, we have 

\begin{eqnarray*}
	|f_h (x, {\bf y})| \leqslant & \displaystyle \frac{\norm{\frac{\partial h}{\partial x}}{\infty}}{ - \int_x^{- \nu} u q_{\nu} (u) \ \mathrm{d}u} & \left[  \int_{x}^{- \nu} \left[\tilde{F}_{\nu} (x) - \tilde{F}_{\nu} (w) \right] \ \mathrm{d}w + \tilde{F}_{\nu}(x)  \int_{- \nu}^{+ \infty} (1 - F_{\nu} (w)) \ \mathrm{d}w \right] \\
    \leqslant & \displaystyle \frac{\norm{\frac{\partial h}{\partial x}}{\infty}}{ - \int_x^{- \nu} u q_{\nu} (u) \ \mathrm{d}u} & \left[ \tilde{F}_{\nu} (x)  \left( - (x + \nu) + \int_{- \nu}^{+ \infty} (1 - F_{\nu} (w)) \ \mathrm{d}w \right) - \int_x^{- \nu} \tilde{F}_{\nu} (w) \ \textrm{d}w \right] \\
    \leqslant & \displaystyle \frac{\norm{\frac{\partial h}{\partial x}}{\infty}}{ - \int_x^{- \nu} u q_{\nu} (u) \ \mathrm{d}u} & \left[ -x \tilde{F}_{\nu} (x) -  \left( \int_x^{- \nu} u q_{\nu} (u) \ \textrm{d}u - x \tilde{F}_{\nu} (x) \right) \right] \\
    = & \displaystyle \norm{\frac{\partial h}{\partial x}}{\infty}. & 
\end{eqnarray*}

We used Lemma \ref{ll3} to express the integral of $\tilde{F}_{\nu}$, and the fact that $\int_{- \nu}^{+ \infty} (1 - F_{\nu} (u)) \ \textrm{d}u = \nu$. Hence, we proved the bound (\ref{ingfh}) on $(- \infty, - \nu)$. Using the continuity  of $f_h$, this inequality also holds for  $x = - \nu$. So, the inequality (\ref{ingfh}) is true. 

\vskip0.2cm Let us now prove the inequalities (\ref{ingdfh}). We start by dealing with the derivative of $f_{h}$ with respect to $x$. Since $f_{h}$ solves (\ref{s2}), we have for $x \neq - \nu$ and for ${\bf y}\in \mathbb{R} ^{n}$, 
\begin{equation}
	\label{24i-2}
 \frac{\partial f_h}{\partial x} (x, \mathbf{y}) = \frac{x}{2(x+\nu)} f_h (x, \mathbf{y}) + \frac{h(x, \mathbf{y}) - \E[h(Z_{\nu}, \mathbf{y})]}{2 (x + \nu)}. 
\end{equation}
We again separate the cases $x>-\nu$ and $ x<-\nu$. 

\vskip0.2cm

\noindent $\bullet$ \textit{Assume  $x > -\nu$.} By Lemma \ref{ll3}, point 1. and by relation (\ref{new1}) in the proof of this lemma, we get 

\begin{eqnarray*}
&&	x f_h (x, \mathbf{y}) + \tilde{h}(x, \mathbf{y}) \\
& =& \frac{x(1 - F_{\nu} (x))}{\int_{- \nu}^x u p_{\nu} (u) \ \textrm{d}u} \int_{- \nu}^x \frac{\partial h}{\partial x} (w, \mathbf{y}) F_{\nu} (w) \ \textrm{d}w +
	 \frac{x F_{\nu} (x)}{\int_{- \nu}^x u p_{\nu} (u) \ \textrm{d}u} \int_x^{+ \infty} \frac{\partial h}{\partial x} (w, \mathbf{y}) (1 - F_{\nu} (w)) \ \textrm{d}w \\
	& + & \int_{- \nu}^{x} \frac{\partial h}{\partial x} (w, \mathbf{y}) F_{\nu} (w) \ \textrm{d}w - \int_{x}^{+ \infty} \frac{\partial h}{\partial x} (w, \mathbf{y})  (1 - F_{\nu} (w))  \ \textrm{d}w \\
	&= & \frac{G(x)}{A_{\nu} (x)} \int_{- \nu}^x \frac{\partial h}{\partial x} (w, \mathbf{y}) F_{\nu} (w) \ \textrm{d}w - \frac{H(x)}{A_{\nu} (x)} \int_x^{+ \infty} \frac{\partial h}{\partial x} (w, \mathbf{y}) (1 - F_{\nu} (w)) \ \textrm{d}w,
\end{eqnarray*}
\noindent where  we used the notation 
$$ \left\{ \begin{array}{rcl}
	A_{\nu} (x) & = & \int_{- \nu}^x u p_{\nu} (u) \ \textrm{d}u = - 2 (x+\nu) p_{\nu} (x) \leqslant 0  \\
	G(x) & = & A_{\nu} (x) + x(1 - F_{\nu} (x)) \\ 
	H(x) & = & - A_{\nu} (x) + x F_{\nu} (x).
\end{array} \right. $$
\noindent  Notice that $H'(x) = F_{\nu} (x)$ and $H(- \nu) = 0$. Hence, $H(x) = \int_{- \nu}^x F_{\nu} (u) \ \textrm{d}u$. We also have $G(x) = \int_x^{+ \infty} (1 - F_{\nu} (u)) \ \textrm{d}u$, since $G'(x) = 1 - F_{\nu} (x)$, $G(- \nu) = - \nu$ and $\int_{- \nu}^{+ \infty} (1 - F_{\nu}(u)) \ \textrm{d}u = \nu$. Therefore, we have the following bound for  $\frac{\partial f_h}{\partial x}$  :

$$ \left| \frac{\partial f_h}{\partial x} (x, \mathbf{y}) \right| \leqslant \norm{\frac{\partial h}{\partial x}}{\infty} \frac{2 \left( \int_{- \nu}^x F_{\nu} (u) \ \textrm{d}u \right) \left( \int_{x}^{+ \infty} (1-F_{\nu} (u)) \ \textrm{d}u \right)}{(x+\nu) | A_{\nu} (x) |}. $$
To conclude, we need to prove that the function $S_{\nu}$ defined by

\begin{equation}
\label{24i-1} S_{\nu} (x) := \frac{2 \left( \int_{- \nu}^x F_{\nu} (u) \ \textrm{d}u \right) \left( \int_{x}^{+ \infty} (1-F_{\nu} (u)) \ \textrm{d}u \right)}{(x+\nu) |A_{\nu} (x)|} 
\end{equation}
\noindent is bounded on $(- \nu, + \infty)$. We will prove that $S_{\nu}$ is bounded on $\R_+$ and admits a finite limit at $- \nu$. For $x > 0$, we use Lemma \ref{ll3} to have 

\begin{eqnarray*}
    \int_{x}^{+ \infty} (1 - F_{\nu} (u)) \ \textrm{d}u = & & 2(x+\nu) p_{\nu}(x) - x(1 - F_{\nu}(x)) \\
    \leqslant & & 2(x+\nu)p_{\nu} (x).
\end{eqnarray*}

\noindent Consequently, for every $x > 0$, we have (by using that $F_{\nu}$ is increasing and $\int_{-\nu}^{\infty} (1-F(u))du=\nu$), 

$$ S_{\nu} (x) \leqslant \frac{2(x+\nu) F_{\nu} (x) \cdot 2 (x+\nu) p_{\nu} (x)}{2(x+\nu)^2 p_{\nu} (x)} = 2 F_{\nu} (x) \leqslant 2. $$

\noindent Now, for $x$ going to $- \nu$, notice first that :

$$ \lim_{x \to - \nu} S_{\nu} (x) =   2\nu \lim_{x \to - \nu} \frac{\int_{-\nu}^x F_{\nu} (u) \ \textrm{d}u}{2(x+\nu)p_{\nu}(x) }, $$

\noindent Since both numerator and denominator go to zero, we apply l'Hôpital's rule several times to conclude. We use the  relation (\ref{5mm-1}) to differentiate in (\ref{24i-1}). 
\begin{eqnarray*}
    \lim_{x \to - \nu} S_{\nu} (x) = & & 2 \nu \lim_{x \to - \nu} \frac{F_{\nu} (x)}{(2-x)(x+\nu)p_{\nu}(x)} = \frac{4 \nu }{2 + \nu} \lim_{x \to - \nu} \frac{F_{\nu} (x)}{2(x+\nu) p_{\nu} (x)} \\
    = & &  \frac{4 \nu }{2 + \nu} \lim_{x \to - \nu} \frac{p_{\nu} (x)}{-x p_{\nu}(x)} = \frac{4}{2+\nu}.
\end{eqnarray*}

\noindent Hence, $S_{\nu}$ is continuous, bounded on $\R_+$ and admits a finite limit at $-\nu$, so $S_{\nu}$ is bounded in $(- \nu, + \infty)$ and so we conclude that  

$$ \forall x > - \nu,  \left| \frac{\partial f_h}{\partial x} (x, \mathbf{y}) \right| \leqslant \norm{S_{\nu}}{\infty} \norm{\frac{\partial h}{\partial x}}{\infty}. $$

\vskip0.2cm

\noindent $\bullet$ \textit{For $x < - \nu$.} We follow the same computations. We get  
\begin{eqnarray*}
	& & x f_{h} (x, \mathbf{y}) + \tilde{h}(x, \mathbf{y}) \\
	= & & \frac{\tilde{G}(x)}{B_{\nu} (x)} \int_x^{+ \infty} \frac{\partial h}{\partial x} (w, \mathbf{y}) \left[ 1 - F_{\nu} (w) \right] \ \textrm{d}w -  \frac{x}{B_{\nu} (x)} \int_{x}^{- \nu} \frac{\partial h}{\partial x} (w, \mathbf{y}) \tilde{F}_{\nu} (w) \ \textrm{d}w, 
\end{eqnarray*}

\noindent where we denoted

$$ \left\{ \begin{array}{rcl}
	B_{\nu} (x) & = & \int_{x}^{- \nu} u q_{\nu} (u) \ \textrm{d}u = 2 (x+\nu) q_{\nu} (x)  \\
	\tilde{G}(x) & = & x \tilde{F}_{\nu} (x) - B_{\nu} (x). \\
\end{array} \right. $$

\noindent Then, we have $\tilde{G}(x) = \int_x^{- \nu} \tilde{F}_{\nu} (u) \ \textrm{d}u$.  Set for all $x < - \nu$ :

$$ R_{\nu} (x) := \frac{\left( \int_x^{- \nu} \tilde{F}_{\nu} (u) \ \textrm{d}u \right) \left( \int_x^{+ \infty} \left[ 1 - F_{\nu} (w) \right] \ \textrm{d}w - x \right)}{(x+ \nu) B_{\nu} (x)} = \frac{-2x \left( \int_x^{- \nu} \tilde{F}_{\nu} (u) \ \textrm{d}u \right) }{2 q_{\nu} (x) (x+ \nu)^2} .$$
Then we have by (\ref{24i-2}), for every $x<-\nu$ and for every $\mathbf{y} \in \mathbb{R}^{n}$, 

$$  \left| \frac{\partial f_h}{\partial x} (x, \mathbf{y}) \right| \leqslant R_{\nu} (x) \norm{\frac{\partial h}{\partial x}}{\infty} . $$

We need to prove that $R_{\nu}$ is bounded. Here, we will prove that $R_{\nu}$ admits finite limits at $- \infty$ and at $-\nu$. We will use  again the fundamental equality  (\ref{5mm-1}).  In $- \nu$ :
\begin{eqnarray*}
    \lim_{x \to - \nu} R_{\nu} (x) = & &  2 \nu \lim_{x \to - \nu} \frac{\int_x^{- \nu} \tilde{F}_{\nu} (u) \ \textrm{d}u}{2(x+\nu)^2 q_{\nu} (x)} \\
    = & & 2 \nu \lim_{x \to - \nu} \frac{- \tilde{F}_{\nu} (x) }{(2-x)(x+\nu) q_{\nu} (x)} = \frac{4 \nu}{\nu + 2} \lim_{x \to - \nu} \frac{- \tilde{F}_{\nu} (x) }{2(x+\nu) q_{\nu} (x)} \\
    = & & \frac{4 \nu}{\nu + 2} \lim_{x \to - \nu} \frac{q_{\nu} (x)}{-x q_{\nu} (x)} = \frac{4}{2+\nu}.
\end{eqnarray*}

\noindent In $- \infty$, we proceed in a similar way, 
\begin{eqnarray*}
    \lim_{x \to - \infty} R_{\nu} (x) = & & -2 \lim_{x \to - \infty} \frac{\int_x^{- \nu} \tilde{F}_{\nu} (u) \ \textrm{d}u}{2(x+\nu) q_{\nu}(x)} \\
    = & &  -2 \lim_{x \to - \infty} \frac{ \tilde{F}_{\nu} (x) }{x q_{\nu}(x)} = -4 \lim_{x \to - \infty} \frac{ \tilde{F}_{\nu} (x) }{2 (x+ \nu) q_{\nu}(x)} \\
    = & & -4 \lim_{x \to - \infty} \frac{q_{\nu} (x)}{x q_{\nu}(x)} = 0.
\end{eqnarray*}

\noindent This proves that $R_{\nu}$ is bounded on $(- \infty, - \nu)$, and so that for every $x < -\nu, {\bf y} \in \mathbb{R} ^{n}$,

$$ \left| \frac{\partial f_h}{\partial x} (x, \mathbf{y}) \right| \leqslant \norm{R_{\nu}}{\infty} \norm{\frac{\partial h}{\partial x}}{\infty}. $$

\noindent Since $f_h$ is $C^1$ with respect to $x$, this holds on $- \nu$ by continuity. Hence, we proved that 

$$ \norm{\frac{\partial f_h}{\partial x}}{\infty} \leqslant \max \left\{ \norm{S_{\nu}}{\infty}, \norm{R_{\nu}}{\infty} \right\} \norm{\frac{\partial h}{\partial x}}{\infty}$$
and the first inequality in (\ref{ingdfh}) is done. 

To deal with the second bound in (\ref{ingdfh}), by differentiating the Stein's equation (\ref{s2}) with respect to $y_j$, we observe that we have in fact 

$$ \frac{\partial f_h}{\partial y_j} = f_{\frac{\partial h}{\partial y_j}}.  $$

\noindent Since $h$ is $C^2$, the estimation made in (\ref{ingfh}) directly yields  the second estimation of (\ref{ingdfh}).
\qed

\section{Multidimensional Stein-Malliavin calculus for the Gamma distribution}\label{sec3}

We use the multidimensional Stein's equation and the properties of its solution proven in Proposition \ref{pp2} in combination with the techniques of Malliavin calculus in order to  obtain bounds for the distance between the probability distributions of random vectors. Actually, we will evaluate the second Wasserstein distance between the law of an arbitrary random vector $(X,\mathbb{Y})$ in $\mathbb{R}\times \mathbb{R}^{n}$ and the random vector $ (Z_{\nu}, \mathbb{Y})$, where $ Z_{\nu}\sim F(\nu)$ and $Z_{\nu}$ is independent by $\mathbb{Y}$. 

Let us first introduce the distances that we use in the sequel. Let $\mathcal{H}$ the set of functions $h : \R^{n+1}  \longrightarrow \R$ which are $C^1$ such that every partial derivative of $h$ are Lipschitz with 

$$ \norm{h}{\text{Lip}} + \norm{h'}{\text{Lip}} \leqslant 1, $$

\noindent where 

$$ \norm{h}{\text{Lip}} := \sup_{\substack{{\bf x,y}\in \R^{n+1} \\ \mathbf{x} \neq \mathbf{y}}} \frac{|h(\mathbf{x}) - h(\mathbf{y})|}{\norm{\mathbf{x} - \mathbf{y}}{}} \text{ and } \norm{h'}{\text{Lip}} := \max_{1 \leqslant j \leqslant n+1} \sup_{\substack{{\bf x,y}\in \R ^{n+1} \\ \mathbf{x} \neq \mathbf{y}}} \frac{|\partial_j h (\mathbf{x}) - \partial_j h (\mathbf{y})|}{\norm{\mathbf{x} - \mathbf{y}}{}}. $$

Then, for every integrable random vectors $\mathbb{X}$ and $\mathbb{Y}$ of $\R^{n+1}$, we define the \textit{second Wasserstein distance} defined by 

\begin{equation}
\label{d2}
d_2 (\P_{\mathbb{X}}, \P_{\mathbb{Y}}) := \sup_{h \in \mathcal{H}} \left| \E[h(\mathbb{X})] - \E[h(\mathbb{Y})] \right|. 
\end{equation}

\noindent (By convention, if $\mathbb{X}$ or $\mathbb{Y}$ are not integrable, we set it equal to $\infty$: we cannot measure anything with this distance). We also recall the classical Wasserstein distance. Let

$$ \mathcal{A} = \{ h: \mathbb{R}^{n}\to \mathbb{R}, \  h \mbox{ is Lipschitz continuous with } \Vert h\Vert _{\text{Lip}}\leqslant 1 \}.$$
Then the Wasserstein distance between the probability distributions of $\mathbb{X}$ and $\mathbb{Y}$ is defined by 
\begin{equation}
	\label{dw}
	d_{\textrm{W}}(\P_{\mathbb{X}}, \P_{\mathbb{Y}}) : = \sup_{h \in \mathcal{A}} \left| \mathbf{E}[h(\mathbb{X})] - \mathbf{E} [h(\mathbb{Y})] \right|.
\end{equation}
We can show the following fact: the convergence for the distance $d_{2}$  implies the convergence for the classical Wasserstein distance.

\begin{lemma}
	For every integrable random vectors $\mathbb{X}$ and $\mathbb{Y}$, we have 
	
	$$ d_{\mathrm{W}} (\mathbb{X}, \mathbb{Y}) \leqslant \sqrt{\frac{32 I_n}{\sqrt{\pi}}} \sqrt{d_2 (\mathbb{X}, \mathbb{Y})}, $$
	
	\noindent where $I_n := \frac{1}{\sqrt{(2 \pi)^n}} \int_{\R^n} \norm{\mathbf{u}}{} e^{\frac{- \norm{\mathbf{u}}{}^2}{2}} \ \mathrm{d} \mathbf{u}$.
\end{lemma}
\noindent {\bf Proof: } Let $h$ be a $1$-Lipschitzian map on $\R^n$. Then by Rademacher's theorem, $h$ is Lebesgue almost-everywhere differentiable. We can consider a version of its partial derivatives. For $\sigma > 0$, we consider the $n$-dimensional Gaussian  kernel with variance $\sigma ^{2}$ given by 

$$ G_{\sigma} (\mathbf{x}) := \frac{1}{\sqrt{(2 \pi \sigma^2)^n}} \exp \left( \frac{- \norm{\mathbf{x}}{}^2}{2 \sigma^2} \right), $$

\noindent and the following approximation of $h$: $h_{\sigma} := h * G_{\sigma}$.  In a first step,  let us prove that there exists a constant $a_{\sigma}$ such that $a_{\sigma} h_{\sigma} \in \mathcal{H}$. The function $h_{\sigma}$ is $C^{\infty}$ and its first partial  derivatives are given by $\partial_i h_{\sigma} = \partial_i h * G_{\sigma}$. Hence, we conclude that 

$$ \norm{\partial_j h_{\sigma}}{\infty} \leqslant \norm{h}{\text{Lip}} \leqslant 1. $$

\noindent Concerning  the second derivative of $h_{\sigma}$, we have:

\begin{eqnarray*}
	\left| \partial^2_{i,j} h_{\sigma} (\mathbf{x}) \right| = & & \left| \partial_i h * \partial_j G_{\sigma} (\mathbf{x}) \right| \leq   \frac{\norm{h}{\text{Lip}}}{\sigma^2} \int_{\R^n} |u_j| \exp \left( \frac{- \norm{\mathbf{u}}{}^2}{2 \sigma^2} \right) \frac{\text{d}\mathbf{u}}{\sqrt{(2 \pi \sigma^2)^n}} \\
	 &= & \sqrt{\frac{2}{\pi}} \frac{\norm{h}{\text{Lip}}}{\sigma} \leq  \sqrt{\frac{2}{\pi \sigma^2}}.
\end{eqnarray*}
\noindent  This implies that  for every $\sigma \in \left[ 0, \sqrt{\frac{2}{\pi}} \right]$, $\sqrt{\frac{\pi \sigma^2}{2}} h_{\sigma} \in \mathcal{H}$.

 Let us  next evaluate the difference between  $h$ and $h_{\sigma}$. We can write

\begin{eqnarray*}
&&	\left| h(\mathbf{x}) - h_{\sigma} (\mathbf{x}) \right| \\
 &= & \left| \int_{\R^n} (h(\mathbf{x}) - h(\mathbf{y})) G_{\sigma} (\mathbf{x} - \mathbf{y}) \frac{\text{d}\mathbf{y}}{\sqrt{(2 \pi \sigma^2)^n}} \right| \leq  \norm{h}{\text{Lip}} \int_{\R^n } \norm{\mathbf{u}}{} G_{\sigma} (\mathbf{u}) \frac{\text{d}\mathbf{u}}{\sqrt{(2 \pi \sigma^2)^n}} \\ 
 &\leq  & \frac{\sigma^{n+1} \norm{h}{\text{Lip}}}{\sqrt{(2 \pi \sigma^2)^n}} \int_{\R^n} \norm{\mathbf{u}}{} e^{\frac{- \norm{\mathbf{u}}{}^2}{2}} \ \textrm{d} \mathbf{u} \leq  \left( \frac{1}{\sqrt{(2 \pi)^n}} \int_{\R^n} \norm{\mathbf{u}}{} e^{\frac{- \norm{\mathbf{u}}{}^2}{2}} \ \textrm{d} \mathbf{u} \right) \sigma := I_n \sigma.
\end{eqnarray*}
Therefore, for every $\sigma \in \left[ 0, \sqrt{\frac{2}{\pi}} \right]$, we have 

\begin{eqnarray*}
	\left| \E[h(\mathbb{X})] - \E [h(\mathbb{Y})] \right| \leqslant & & \left| \E[h(\mathbb{X})] - \E [h_{\sigma} (\mathbb{X})] \right| + \left| \E[h(\mathbb{Y})] - \E [h_{\sigma} (\mathbb{Y})] \right| \\
	& + & \left| \E[h_{\sigma} (\mathbb{X})] - \E [h_{\sigma} (\mathbb{Y})] \right| \\
	\leqslant & & 2 \norm{h - h_{\sigma}}{\infty} + \sigma \sqrt{\frac{\pi}{2}} d_2 (X, Y) \\
	\leqslant & &  2 I_n \sigma + \frac{1}{\sigma} \sqrt{\frac{2}{\pi}} d_2 (\mathbb{X}, \mathbb{Y}).
\end{eqnarray*}

\noindent The left hand side, seen as a function of $\sigma$, can be optimized on $\left[ 0, \sqrt{\frac{\pi}{2}} \right]$ by taking 

$$ \sigma_0 = \sqrt{\frac{\sqrt{\frac{2}{\pi}} d_2 (\mathbb{X}, \mathbb{Y})}{2 I_n} } \in \left[ 0, \sqrt{\frac{\pi}{2}} \right]. $$

\noindent Note that since $(I_n)_n$ is increasing, we have indeed $\sigma_0 \leqslant \sqrt{\frac{d_2 (\mathbb{X}, \mathbb{Y})}{2}} \leqslant \sqrt{\frac{\pi}{2}}$, since the distance is lower than $1$. This value of $\sigma$ yields to  

$$ d_{\textrm{W}} (\mathbb{X}, \mathbb{Y}) \leqslant 2 \sqrt{\frac{8 I_n}{\sqrt{\pi}} d_2 (\mathbb{X}, \mathbb{Y})}, $$
which is the desired conclusion. \qed 

Let us state and prove the main result of this section. 

\begin{theorem}
\label{tt1}	
	Let $X$ a centered random variable and $\mathbb{Y} = (Y_1, \cdots, Y_n)$ a centered random vector of $\R^n$ such that $X \in \mathbb{D}^{1, 2}$ and $Y_j \in \mathbb{D}^{1, 2}$ for every $j \in \{ 1, \cdots ,n \}$. Then there exists a constant $C = C (\nu) > 0$ such that
	
	\begin{eqnarray}{\label{distL1}}
		& & d_{2} \left( \P_{(X, \mathbb{Y})}, F (\nu) \otimes \P_{\mathbb{Y}} \right) \\ 
	 &\leq  & C \ \E \left[ \left| 2(X + \nu) - \left\langle \mathrm{D}(-L)^{-1} X, \mathrm{D}X  \right\rangle \right| \right] + C \ \sum_{j=1}^n \E \left[ \left| \left\langle \mathrm{D} (-L)^{-1} X, \mathrm{D} Y_j \right\rangle \right| \right]. \nonumber
	\end{eqnarray}
	
	\noindent If moreover we suppose that $X \in \mathbb{D}^{1, 4}$ and $Y_j \in \mathbb{D}^{1, 4}$ for every $j \in \{ 1, \cdots ,n \}$,  then
	
	\begin{eqnarray}{\label{distL2}}
		& & d_{2} \left( \P_{(X, \mathbb{Y})}, F (\nu) \otimes \P_{\mathbb{Y}} \right) \\ 
	&\leq  & C \ \E \left[ \left( 2(X + \nu) - \left\langle \mathrm{D}(-L)^{-1} X, \mathrm{D}X  \right\rangle \right)^2 \right]^{\frac{1}{2}} + C \ \sum_{j=1}^n \E \left[ \left\langle \mathrm{D} (-L)^{-1} X, \mathrm{D} Y_j \right\rangle^2 \right]^{\frac{1}{2}}. \nonumber
	\end{eqnarray}
	
\end{theorem}
\noindent{\bf Proof: }  Let suppose first that $h \in \mathcal{C}^2 (\R^n \times \R)$ with $\norm{h}{\text{Lip}} + \norm{h'}{\text{Lip}} \leqslant 1$. Let  $f_h$ be the corresponding solution to the Stein's equation (\ref{s2}).  Then 

\begin{eqnarray*}
    & & \int_{(- \nu, + \infty) \times \R^n} h(x, \mathbf{y}) \ \textrm{d}\P_{(X, \mathbb{Y})} (x, \mathbf{y}) - \int_{(- \nu, + \infty) \times \R^n} h(x, \mathbf{y}) \ \textrm{d}\P_{Z_{\nu}} \otimes \textrm{d}\P_{\mathbb{Y}} (x, \mathbf{y}) \\
    = & & \int_{(- \nu, + \infty) \times \R^n} h(x, \mathbf{y}) \ \textrm{d}\P_{(X, \mathbb{Y})} (x, \mathbf{y}) - \int_{\R^n} \E[h(Z_{\nu}, \mathbf{y})] \ \textrm{d}\P_{\mathbb{Y}} (\mathbf{y}) \\
    = & & \E \left[ h(X, \mathbb{Y}) - \E \left[ h(Z_{\nu}, \mathbb{Y}) \right] \right].
\end{eqnarray*}

So, by the Stein's equation (\ref{s2}), we have 

\begin{eqnarray}{\label{thm1}}
	& & \int_{(- \nu, + \infty) \times \R^n} h(x, \mathbf{y}) \ \textrm{d}\P_{(X, \mathbb{Y})} (x, \mathbf{y}) - \int_{(- \nu, + \infty) \times \R^n} h(x, \mathbf{y}) \ \textrm{d}\P_{Z_{\nu}} \otimes \textrm{d}\P_{\mathbb{Y}} (x, \mathbf{y}) \nonumber \\
	= & & 2 \E \left[ (X + \nu) \frac{\partial f_h}{\partial x} (X, \mathbb{Y}) \right] - \E \left[ X f_h (X, \mathbb{Y}) \right]. 
\end{eqnarray}
Since $X$ is centered, we can write $ X = \delta \left[ \textrm{D} (-L)^{-1} X \right]$, and by plugging this  into (\ref{thm1})

\begin{eqnarray*}
	& & \int_{(- \nu, + \infty) \times \R^n} h(x, \mathbf{y}) \ \textrm{d}\P_{(X, \mathbb{Y})} (x, \mathbf{y}) - \int_{(- \nu, + \infty) \times \R^n} h(x, \mathbf{y}) \ \textrm{d}\P_{Z_{\nu}} \otimes \textrm{d}\P_{\mathbb{Y}} (x, \mathbf{y}) \\
	= & & 2 \E \left[ (X + \nu) \frac{\partial f_h}{\partial x} (X, \mathbb{Y}) \right] - \E \left[ \delta \left[ \textrm{D} (-L)^{-1} X \right] f_h (X, \mathbb{Y}) \right] \\
	= & & 2 \E \left[ (X + \nu) \frac{\partial f_h}{\partial x} (X, \mathbb{Y}) \right] - \E \left[  \left\langle \textrm{D} (-L)^{-1} X , \textrm{D} f_h (X, \mathbb{Y}) \right\rangle \right] \\
	= & &  2 \E \left[ (X + \nu) \frac{\partial f_h}{\partial x} (X, \mathbb{Y}) \right] - \E \left[ \frac{\partial f_h}{\partial x} (X, \mathbb{Y}) \left\langle \textrm{D} (-L)^{-1} X , \textrm{D} X \right\rangle \right] \\
	& - &  \sum_{j=1}^n \E \left[ \frac{\partial f_h}{\partial y_j} (X, \mathbb{Y}) \left\langle \textrm{D} (-L)^{-1} X , \textrm{D} Y_j \right\rangle \right].
\end{eqnarray*}
Since $\norm{h}{\text{Lip}} + \norm{h'}{\text{Lip}} \leqslant 1$, we have $\norm{\frac{\partial f_h}{\partial x}}{\infty} \leqslant C$ and $\norm{\frac{\partial f_h}{\partial y_j}}{\infty} \leqslant 1$ by using Lemma \ref{ll3}.  So, we obtain for a generic constant $C$ 

\begin{eqnarray*}
	& & \left| \int_{(- \nu, + \infty) \times \R^n} h(x, \mathbf{y}) \ \textrm{d}\P_{(X, \mathbb{Y})} (x, \mathbf{y}) - \int_{(- \nu, + \infty) \times \R^n} h(x, \mathbf{y}) \ \textrm{d}\P_X \otimes \textrm{d}\P_{\mathbb{Y}} (x, \mathbf{y}) \right| \\
	\leqslant & & C \E \left[ \left| (X + \nu) - \left\langle \textrm{D} (-L)^{-1} X , \textrm{D} X \right\rangle \right| \right] + C \sum_{j=1}^n \E \left[ \left| \left\langle \textrm{D} (-L)^{-1} X , \textrm{D} Y_j \right\rangle \right| \right].   \\
\end{eqnarray*}
To conclude the bound  (\ref{distL1}), we need to have the above  inequality for every $h \in \mathcal{H}$. If $h$ is one of those functions, we approach it by the sequence

$$ h_k (x, \mathbf{y}) := \E \left[ h \left( x + \frac{N}{\sqrt{k}}, \mathbf{y} + \frac{\mathbf{N}}{\sqrt{k}} \right) \right], $$
where $N \sim \mathcal{N}(0, 1)$ is independent of $\mathbf{N} \sim \mathcal{N}_n (0, I_n)$. Then $(h_k)_k$ uniformly converges to $h$ and we still have $\norm{h_k}{\text{Lip}} + \norm{h_k'}{\text{Lip}} \leqslant 1$. Consequently, we have:

\begin{eqnarray*}
	& & \left| \int_{(- \nu, + \infty) \times \R^n} h(x, \mathbf{y}) \ \textrm{d}\P_{(X, \mathbb{Y})} (x, \mathbf{y}) - \int_{(- \nu, + \infty) \times \R^n} h(x, \mathbf{y}) \ \textrm{d}\P_{Z_{\nu}} \otimes \textrm{d}\P_{\mathbb{Y}} (x, \mathbf{y}) \right| \\
	\leqslant & & 2 \norm{h - h_k}{\infty} + C \E \left[ \left| (X + \nu) - \left\langle \textrm{D} (-L)^{-1} X , \textrm{D} X \right\rangle \right| \right]   \\
	& + & C \sum_{j=1}^n \E \left[ \left| \left\langle \textrm{D} (-L)^{-1} X , \textrm{D} Y_j \right\rangle \right| \right].
\end{eqnarray*}

Taking the limit as $k \to \infty$, (\ref{distL1}) is obtained. For (\ref{distL2}), it suffices to notice that if $X, Y_{j} \in \mathbb{D}^{1, 4}$,  then the  scalar products $\left\langle \textrm{D} (-L)^{-1} X , \textrm{D} X \right\rangle$ and $\left\langle \textrm{D} (-L)^{-1} X, \textrm{D} Y_j \right\rangle$ are in $L^2$, so that we can apply Cauchy-Schwarz's  inequality. \qed

\section{Asymptotic independence on Wiener chaos}\label{sec4}
We now focus on random variables in Wiener chaos and we give an asymptotic variant of the result proven in Theorem \ref{tt1}. We will consider a sequence $(X_{k}, k\geq 1)$ that converges in distribution  when $k\to \infty$  to the centered Gamma law $ F(\nu)$ and a sequence of random vectors $(\mathbb{Y}_{k}, k\geq 1)$ converging  in law as $k\to\infty$ to an arbitrary  random vector $\mathbb{Y}$.  We assume that for each $k\geq 1$,  $X_{k}$ and the components of $ \mathbb{Y}_{k}$ belong to a Wiener chaos of fixed order. Under some pretty natural assumptions, we deduce, by using Theorem \ref{tt1} and the properties of random variables in Wiener chaos, that the random sequence  $\left( (X_{k}, \mathbb{Y}_{k}), k\geq 1\right)$ converges in law to $ (Z_{\nu}, \mathbb{Y})$, where $Z_{\nu}$ follows the centered Gamma distribution with  parameter $\nu>0$ and $Z_{\nu}$ and $\mathbb{Y}$ are independent. This means that the sequences $ (X_{k}, k\geq 1)$ and $(\mathbb{Y}_{k}, k\geq 1)$ are asymptotically independent. We obtain bounds to quantify this asymptotic independence under the $d_{2}$-distance defined by (\ref{d2}). 

Let us first recall, that if a sequence $(X_{k}, k\geq 1)$ in the $q$th  Wiener chaos converges to $F(\nu)$, then the order $q$ of the chaos must be an even integer. Indeed, if $q$ is odd then we have $ \E \left[ X_{k}^{3} \right]=0$ for any $k \geqslant 1$ and it contradicts the fact that $\E \left[ Z_{\nu}^{3} \right] = 8 \nu >0$ if $Z_{\nu}\sim F(\nu )$. 

Before stating our result, let us recall the following criterion for the Gamma approximation on Wiener chaos. We refer to \cite{NP2} for its proof. 

\begin{theorem}\label{teo-np}
	Let $ (X_{k}, k\geq 1)$ be a sequence of random variables such that for every $k \geqslant 1$, $X_{k}=I_{q} ( f_{k})$ where $ q$ is an even integer and $f_{k} \in H ^{\odot q}$. Assume that $\E \left[ X_{k}^{2} \right] \xrightarrow[k \to \infty]{} 2\nu$. Then the following are equivalent:
	\begin{enumerate}
		\item The sequence $(X_{k}, k\geq 1)$ converges in distribution to $F(\nu)$. 
		
		\item For every $p \in \{ 1,..., q-1 \}$ with $p \neq \frac{q}{2}$,
		\begin{equation*}
			\norm{ f_{k}\otimes _{p} f_k}{H ^{\otimes 2q-2p}}^{2} \xrightarrow[k \to \infty]{} 0 \ \text{ and } \ \norm{f_{k} \, \widetilde{ \otimes}_{\frac{q}{2}} \, f_{k} -c_{q} f_{k}}{H ^{\otimes q}}^{2} \xrightarrow[k \to \infty]{} 0,
		\end{equation*} 
where $c_{q}>0$ is an explicit constant. 
\item  As $k\to \infty$, 
\begin{equation*}
	 \norm{\mathrm{D} X_{k}}{H}^2-2qX_{k} -2q\nu \xrightarrow[k \to \infty]{} 0 \mbox{ in } L ^{2} (\Omega).
\end{equation*}
	\end{enumerate}
\end{theorem}   
\noindent Other equivalent conditions to 1.-3. (not needed in our work) are stated and proven in \cite{NP2}.

We apply this theorem in our context, about asymptotic independence with a component going to a centered Gamma law. We have the following result. 

\begin{theorem}{\label{thmchaos}}
	Let $n$, $q$ and $q_1, \cdots, q_n$ strictly positive integers such that $q$ is even and $q \geqslant q_j$, for every $j \in \{ 1, \cdots ,n \}$. We consider a sequence $(X_{k}, k \geqslant 1)= \left( I_{q}(f_{k}), k \geqslant 1\right)$ converging  in distribution, as $k$ goes to $+ \infty$, to a random variable $X$ following $F(\nu)$. We consider a sequence of random vectors $(\mathbb{Y}_k)_k = ((Y^1_k, \cdots, Y^n_k))_k$ such that  $Y_k^j = I_{q_j} (g_k^j)$ and  and there exists a random vector $\mathbb{Y} = (Y^1, \cdots, Y^n)$ such that  $(\mathbb{Y}_k, k \geqslant 1)$ converges in distribution to  $\mathbb{Y}$. We suppose that, for every $j \in \{ 1, \cdots, n \}$ :
	
	\begin{equation}{\label{cdn}}
		\E \left[ X_k Y_k^j \right] \xrightarrow[k \to \infty]{} 0 \text{ and } \norm{f_k \otimes_{\frac{q}{2}} g_k^j}{} \xrightarrow[k \to \infty]{} 0.
	\end{equation}
	Then we have the following convergence in distribution : 
	
	\begin{equation}{\label{cvg}}
		(X_k, \mathbb{Y}_k) \xrightarrow[k \to + \infty]{\mathrm{(d)}} (X, \mathbb{Y}),
	\end{equation}
	
	\noindent where $X \sim F (\nu)$ and $X$ is independent of $\mathbb{Y}$. Moreover, if we denote
	
	$$ \theta_k := \P_{(X_k, \mathbb{Y}_k)} \text{ and } \eta := \P_X \otimes \P_{\mathbb{Y}}, $$
	
	\noindent then we have the following estimation of the second Wasserstein distance between those two measures:
	
	\begin{eqnarray}{\label{ingchaos}}
		d_{2} (\theta_k, \eta) \leqslant & & C \ \E \left[ \left( 2(X_k + \nu) - \left\langle \mathrm{D} (-L)^{-1} X_k, \mathrm{D} X_k \right\rangle \right)^2 \right]^{\frac{1}{2}} \\ 
		& + & C \ \sum_{j=1}^n \E \left[  \left\langle \mathrm{D} (-L)^{-1} X_k, \mathrm{D} Y_k^j \right\rangle^2 \right]^{\frac{1}{2}} + d_{2} \left( \P_{\mathbb{Y}_k}, \P_{\mathbb{Y}} \right). \nonumber 
	\end{eqnarray}
\end{theorem}
\noindent{\bf Proof: } For $k\geq 1$, let  us consider the probability measure $\eta_k := \P_X \otimes \P_{\mathbb{Y}_k} = F (\nu) \otimes \P_{\mathbb{Y}_k}$.  Then we have for every $h\in \mathcal{H}$ (the set defined at the beginning of Section \ref{sec3}), 

\begin{eqnarray*}
	& & \left| \int_{\R} \int_{\R^n} h(x, \mathbf{y}) \ \textrm{d} \eta_k (x, \mathbf{y}) - \int_{\R} \int_{\R^n} h(x, \mathbf{y}) \ \textrm{d} \eta (x, \mathbf{y}) \right| \\
	= & & \left| \int_{\R^n} \left( \int_{\R} h(x, \mathbf{y}) \ \textrm{d} \P_X (x) \right) \textrm{d}\P_{\mathbb{Y}_k} (\mathbf{y}) - \int_{\R^n} \left( \int_{\R} h(x, \mathbf{y}) \ \textrm{d} \P_X (x) \right) \textrm{d}\P_{\mathbb{Y}} (\mathbf{y}) \right| \\ 
	= & & \left| \E \left[ \int_{\R} h(x, \mathbb{Y}_k) \ \textrm{d} \P_X (x) \right] - \E \left[ \int_{\R} h(x, \mathbb{Y}) \ \textrm{d} \P_X (x) \right] \right| \\
	\leqslant & & d_2 (\P_{\mathbb{Y}_k}, \P_{\mathbb{Y}}),
\end{eqnarray*}
the last inequality being true because the function $\mathbf{y} \longmapsto \int_{\mathbb{R} } h(x,\mathbf{y}) \ \textrm{d}x$ also belongs to $\mathcal{H}$ (on $\mathbb{R}^{n}$). Hence,
 \begin{equation}\label{25i-3}
 d_2 (\eta_k, \eta) \leqslant d_2 (\P_{\mathbb{Y}_k}, \P_{\mathbb{Y}}).
\end{equation}
Now by using the triangle inequality, (\ref{25i-3}) and the bound (\ref{distL2}) in Theorem \ref{tt1}, 
\begin{eqnarray*}
	d_{2} (\theta _{k}, \eta) &\leq & d_{2} (\theta _{k}, \eta _{k})+ d_{2} (\eta _{k}, \eta)\\&\leq & C \ \E \left[ \left( 2(X_k + \nu) - \left\langle \mathrm{D} (-L)^{-1} X_k, \mathrm{D} X_k \right\rangle \right)^2 \right]^{\frac{1}{2}} \\ 
	& + & C \ \sum_{j=1}^n \E \left[  \left\langle \mathrm{D} (-L)^{-1} X_k, \mathrm{D} Y_k^j \right\rangle^2 \right]^{\frac{1}{2}} + d_{2} \left( \P_{\mathbb{Y}_k}, \P_{\mathbb{Y}} \right)
\end{eqnarray*}
which is nothing else but (\ref{ingchaos}). Let us now prove that the right-hand side of (\ref{ingchaos}) converges to zero as $ k\to \infty$.  For the first term, the convergence to zero comes from Theorem \ref{teo-np}. Indeed,

\begin{eqnarray*}
	& & \E \left[ \left( 2(X_k + \nu) - \left\langle \mathrm{D} (-L)^{-1} X_k, \mathrm{D} X_k \right\rangle \right)^2 \right] = \frac{1}{q^2} \ \E \left[ \left( \norm{\textrm{D} X_k}{}^2 - 2q X_k - 2q \nu \right)^2 \right] \xrightarrow[k \to \infty]{} 0
\end{eqnarray*}
due to point 3. in Theorem \ref{teo-np}. About the second term, we will use some classical computations about Wiener chaoses. We have:

$$ \E \left[  \left\langle \mathrm{D} (-L)^{-1} X_k, \mathrm{D} Y_k^j \right\rangle^2 \right] = \frac{1}{q^2} \E \left[ \left\langle \textrm{D} X_k, \textrm{D} Y_k^j \right\rangle^2 \right]. $$

\noindent Then, by using the product formula (\ref{prod}), and the assumption  $q \geqslant q_j$, we conclude that 

\begin{eqnarray}{\label{EDXDY}}
	& & \E \left[  \left\langle \mathrm{D} (-L)^{-1} X_k, \mathrm{D} Y_k^j \right\rangle^2 \right] \nonumber \\
	 &= & \frac{1}{q^2} \sum_{r=0}^{q_j - 1} (r!)^2 \binom{q_j - 1}{r}^2 \binom{q-1}{r}^2 (q+q_j - 2 - 2r)! \norm{f_k \ \widetilde{\otimes}_{r+1} \ g_k^j}{}^2 \nonumber \\
& =& \frac{1}{q^2} \sum_{r=1}^{q_j} ((r-1)!)^2 \binom{q_j - 1}{r-1}^2 \binom{q-1}{r-1}^2 (q+q_j - 2r)! \norm{f_k \ \widetilde{\otimes}_{r} \ g_k^j}{}^2.
\end{eqnarray}
\noindent Suppose first that $q > q_j$. We estimate the norm of the symmetrical tensor product.

\begin{eqnarray*}
	\norm{f_k \ \tilde{\otimes}_r \ g_k^j}{}^2 \leqslant & & \norm{f_k \otimes_r g_k^j}{}^2 = \left\langle f_k \otimes_r g_k^j, f_k \otimes_r g_k^j \right\rangle \\
	= & & \left\langle f_k \otimes_{q-r} f_k, g_k^j \otimes_{q_j - r} g_k^j \right\rangle \leq  \norm{f_k \otimes_{q-r} f_k}{} \cdot \norm{g_k^j \otimes_{q_j - r} g_k^j}{} \\
	\leqslant & & \norm{f_k \otimes_{q-r} f_k}{} \cdot \norm{g_j}{}^2 \leq  \frac{\E \left[ \left( Y_k^j \right)^2 \right]}{q_j !} \norm{f_k \otimes_{q-r} f_k}{}.
\end{eqnarray*}

\noindent Hence, we obtain :

\begin{eqnarray*}
	& & \E \left[  \left\langle \mathrm{D} (-L)^{-1} X_k, \mathrm{D} Y_k^j \right\rangle^2 \right] \\
	\leqslant & & \frac{\E \left[ \left( Y_k^j \right)^2 \right]}{q^2 q_j !} \sum_{r=1}^{q_j } ((r-1)!)^2 \binom{q_j - 1}{r-1}^2 \binom{q-1}{r-1}^2 (q + q_j - 2r)! \norm{f_k \otimes_{q-r} f_k}{}.
\end{eqnarray*}

\noindent Then, if for every $j$, $q > q_j$, then by point 2. in Theorem \ref{teo-np} and the second condition of (\ref{cdn}), the whole expression converges to zero as $k\to \infty$ , and so we get  (\ref{cvg}). 

Suppose now that $q = q_j$. We isolate the term $r = q_j = q$ in \ref{EDXDY}, giving : 

\begin{eqnarray*}
	& & \E \left[  \left\langle \mathrm{D} (-L)^{-1} X_k, \mathrm{D} Y_k^j \right\rangle^2 \right] \\
	= & & \frac{1}{q^2} \sum_{r=1}^{q-1 } ((r-1)!)^2 \binom{q - 1}{r-1}^4 (2(q - r))! \norm{f_k \ \widetilde{\otimes}_{r} \ g_k^j}{}^2 + ((q-1)!)^2 \norm{f_k \ \widetilde{\otimes}_q \ g_k^j}{}^2 \\
	= & & \frac{1}{q^2} \E \left[ X_k Y_k^j \right]^2 + \frac{1}{q^2} \sum_{r=1}^{q-1 } ((r-1)!)^2 \binom{q - 1}{r-1}^4 (2(q - r))! \norm{f_k \ \widetilde{\otimes}_{r} \ g_k^j}{}^2.
\end{eqnarray*}

\noindent By both conditions in (\ref{cdn}), and by point 2.  of Theorem \ref{teo-np}, we also conclude in this case that this expectation goes to zero, and so we have (\ref{cvg}). \qed

\begin{remark}
The assumption (\ref{cdn}) is necessary to get the result in Theorem \ref{thmchaos} and it cannot be avoided. To argue this, consider the following trivial example.  Let $ (W(h), h\in H)$ be an isonormal process and let $ h_{1}, h_{2} \in H$ with $\Vert h_{1}\Vert_{H}  =\Vert h_{2}\Vert_{H} =1$ and $\langle h_{1}, h_{2} \rangle _{H}=\rho \in (0,1)$. Define
\begin{equation*}
	X= W(h_{1}) ^{2}-1 = I_{2} ( h_{1} ^{\otimes 2})\mbox{ and } Y = W(h_{2}).
\end{equation*}
Then $X\sim F(1), Y\sim \mathcal{N} (0,1)$ and all the assumptions in the statement of Theorem \ref{thmchaos}, except (\ref{cdn}), are satisfied.  Notice that $ h_{1} ^{\otimes 2} \otimes _{1} h_{2}= \rho h_{1}$ and $ \Vert h_{1} ^{\otimes 2} \otimes _{1} h_{2}\Vert _{H}= \rho \Vert h_{1}\Vert_{H} = \rho \not=0$ so (\ref{cdn}) does not hold true. On the other hand, the components of the vector $(X,Y)$ are  not independent since $h_{1} ^{\otimes 2} \otimes _{1} h_{2}$ does not vanish almost everywhere on $H$, see \cite{UZ}.
\end{remark}

\section{Examples}
We illustrate the results stated in Section \ref{sec4} by two examples. In these examples we consider a two -dimensional sequence of random variables which has on its first component the sequence $U_{n}$ defined below by (\ref{un}) and which converges in law to $ F(1)$. On the second component, we first consider a  fixed random variables in the second Wiener chaos and then another sequence, correlated with $ U_{n}$, which also converges to the centered Gamma law $F(1)$. We obtain the joint convergence of the two-dimensional sequence to a vector with independent components and we derive the associate rate of convergence under the $ d_{2}$-distance.

\subsection{Example 1} Let $ (h_{i}, i\geq 1)$ be orthonormal elements of the Hilbert space $H$. For $n\geq 2$, set 
\begin{equation}\label{un}
	U_{n}= I_{2} \left( \frac{2}{n-1} \sum_{1\leq i<j\leq n} h_{i} \widetilde{\otimes} h{j} \right)=: I_{2} (f_{n}),
\end{equation}
where $I_{2}$ is the multiple integral with respect to an isonormal process $(W(h), h\in H)$. Then the sequence $ (U_{n}, n\geq 1)$ converges in distribution to the centered Gamma law $F(1)$ (see e.g. \cite{AAPS}).  It follows from Section 4 in \cite{AET} that for $n$ large enough,

\begin{equation}
	\label{30i-1}
	\E \left[ \left( 2(U_{n}+1) -\langle \textrm{D} (-L) ^{-1} U_{n}, \textrm{D} U_{n}\rangle \right)^2 \right] \leqslant  \frac{C}{n}. 
\end{equation}
In particular,
\begin{equation}
	\label{30i-2}
	d_{2}(U_{n}, F(1))\leqslant  \frac{C}{\sqrt{n}}.
\end{equation}

\noindent We regard the asymptotic behavior of the two-dimensional sequence $ (U_{n}, G)$ where 

$$ G := 2 H_{2} (W(h_{1}))= I_{2} \left( h_{1}^{\otimes 2} \right). $$ 

\noindent Let us check (\ref{cdn}). Concerning the first part of (\ref{cdn}), we have for every $n\geq 1$,
\begin{eqnarray*}
	\E \left[ U_{n} G \right] &=& \frac{4}{n-1} \sum_{1\leq i<j\leq n} \left\langle h_{i} \widetilde{\otimes} h_{j}, h_{1} ^{\otimes 2} \right\rangle_{H^{\otimes 2}} \\
	&=&\frac{2}{n-1} \sum_{1\leq i<j\leq n} \left\langle h_{i} \otimes  h_{j}+ h_{j}\otimes h_{i}, h_{1} ^{\otimes 2} \right\rangle _{H^{\otimes 2}}\\
	&=&\frac{4}{n-1} \sum_{1\leq i<j\leq n} \mathbf{1}_{\{ i=j=1 \}} =0.
\end{eqnarray*}
Also, 
\begin{eqnarray}
	f_{n} \otimes _{1} h _{1} ^{\otimes 2}&=& \frac{2}{n-1}  \sum_{1\leq i<j\leq n} \left( (h_{i} \widetilde{\otimes} h_{j}) \otimes _{1} h_{1} ^{\otimes 2}\right)\nonumber\\
	&=&  \frac{1}{n-1}  \sum_{1\leq i<j\leq n} \left( (h_{i}\otimes h_{j}+ h_{j}\otimes h_{i})\otimes _{1} h_{1} ^{\otimes 2} \right) \nonumber\\
	&=&  \frac{1}{n-1}  \sum_{1\leq i<j\leq n} \left( (h_{j}\otimes h_{1}) \mathbf{1}_{\{ i=1 \}}+ (h_{i}\otimes h_{1} ) \mathbf{1}_{\{ j=1 \}} \right)=  \frac{1}{n-1}  \sum_{j=2}^n h_{j}\otimes h_{1}\label{30i-4}
\end{eqnarray}
and
\begin{eqnarray}
	\norm{	f_{n} \otimes _{1} h _{1} ^{\otimes 2}}{H ^{\otimes 2}}^{2} &=& \frac{2}{(n-1) ^{2}} \sum_{2 \leqslant j,k \leqslant n} \langle h_{j}\otimes h_{1}, h_{k}\otimes h_{1} \rangle _{H ^{\otimes 2}}\nonumber\\
	&=&  \frac{2}{(n-1) ^{2}} \sum_{2 \leqslant j,k \leqslant n} \langle h_{j}, h_{k} \rangle _{H} =\frac{2(n-1)}{(n-1) ^{2}} \xrightarrow[n \to \infty]{} 0.\label{30i-5}
\end{eqnarray}
It then follows from Theorem \ref{thmchaos} that for $n$ large enough,
\begin{eqnarray}
	\label{30i-3}
	& & d_{2} \left(\P_{(U_{n}, G)}, F(1) \otimes \P_{G} \right) \nonumber \\
	\leqslant & & C \left\{  \E \left[ \left( 2(U_{n}+1) - \left\langle \textrm{D}(-L) ^{-1} U_{n}, \textrm{D}U_{n} \right\rangle \right)^{2} \right]^{\frac{1}{2}} +  \E \left[ \left\langle \textrm{D}(-L) ^{-1} U_{n}, \textrm{D}G \right\rangle^{2} \right]^{\frac{1}{2}} \right\}.
\end{eqnarray}
For any two multiple integrals $X=I_{2}(f)$, $Y=I_{2}(g)$ with $f,g \in H ^{\odot 2}$, 
\begin{equation*}
	\E \left[ \left\langle \textrm{D} (-L)^{-1} X, \textrm{D} Y \right\rangle^{2} \right]= \E [XY]^{2} + 8 \norm{f \otimes_{1} g}{H ^{\otimes 2}}^{2}.
\end{equation*}
Thus, from (\ref{30i-4}) and (\ref{30i-5}), 
\begin{equation}\label{30i-6}
	\E \left[ \left\langle \textrm{D} (-L)^{-1} U_{n}, \textrm{D} G \right\rangle ^{2}\right] = \frac{16}{n-1}.
\end{equation}
By plugging the estimates (\ref{30i-1}) and (\ref{30i-6}) into (\ref{30i-3}), we obtain  for $n$ sufficiently large,
\begin{equation*}
d_{2}\left(\P_{(U_{n}, G)}, F(1) \otimes \P_{G}\right) \leqslant \frac{C}{\sqrt{n}}.
\end{equation*}

\subsection{Example 2}Let $H=L ^{2}(\mathbb{R}_{+})$ and for $i\geq 1$ set
\begin{equation*}
	h_{i} := \mathbf{1}_{[2i, 2i+1]}.
\end{equation*}
Then $(h_{i}, i\geq 1) $ are orthonormal elements in $H$.  Consider the sequence $(U_{n}, n\geq 1)$ given by (\ref{un}). Now define for $i\geq 1$,
\begin{equation*}
	g_{i} := \mathbf{1}_{ \left[ 2i-1+\frac{1}{i}, 2i+\frac{1}{i^a} \right]}
\end{equation*}
with $a >0$. The family $(g_{i}, i\geq 1) $ is also orthogonal in $H$. Notice that for every $i,k \geqslant 1$,
\begin{equation*}
	\label{26i-1}
	\langle h_{i}, g_{k}\rangle _{H} = \frac{1}{i^{a}} \mathbf{1}_{\{ i=k \}}.
\end{equation*}
We consider the sequence $ (V_{n}, n\geq 1)$ given by 
\begin{equation}\label{vn}
	V_{n}= I_{2} \left( \frac{2}{n-1} \sum_{1\leq i<j\leq n} g_{i} \widetilde{\otimes} g_{j} \right)=: I_{2} (\varphi_{n}). 
\end{equation}
We know that both $(U_{n}, n\geq 1)$ (defined at previous example) and $(V_{n}, n\geq 1)$ converge to the centered Gamma law $ F(1)$. We consider the two-dimensional sequence $((U_{n}, V_{n}), n\geq 1)$ and we prove that it converges in law to the vector $(F_{1}, F_{2})$, where $F_{1}, F_{2}\sim F(1)$ and $F_{1}, F_{2}$ are independent. We will also deduce the rate of convergence under the $d_{2}$-distance associated to this limit theorem. 

We compute the quantities in (\ref{cdn}).
\begin{eqnarray}
	\E \left[ U_{n} V_{n} \right] &=& \frac{4}{(n-1)^{2}}\sum_{1\leq i<j\leq n} \sum_{1\leq k<l\leq n} \langle h_{i} \widetilde{\otimes} h_{j},g_{k} \widetilde{\otimes} g_{l}\rangle  \nonumber \\
	&=&  \frac{1}{(n-1) ^{2}}\sum_{1\leq i<j\leq n} \sum_{1\leq k<l\leq n}\langle h_{i}\otimes h_{j}+ h_{j}\otimes h_{i}, g_{k}\otimes g_{l}+ g_{l}\otimes g_{k} \rangle\nonumber \\
	&=& \frac{1}{(n-1) ^{2}}\sum_{1\leq i<j\leq n} \sum_{1\leq k<l\leq n} \frac{2}{i^a j^a} \left( \mathbf{1}_{\{ i=k,j=l \}} + \mathbf{1}_{\{ i=l, j=k \}} \right) \nonumber \\
	&=&\frac{4}{(n-1) ^{2}}\sum_{1\leq i<j\leq n} \frac{1}{(ij)^a} \leq \frac{4}{(n-1) ^{2}}\left( \sum_{i=1}^n \frac{1}{i^a} \right)^{2} \leq \frac{C}{n^{2a}} \xrightarrow[n \to \infty]{} 0. \label{30i-7}
\end{eqnarray}
Next, let us calculate $\norm{ f_{n} \otimes _{1} \varphi_{n}}{ H ^{\otimes 2}}^{2} $. We compute first the contraction :
\begin{eqnarray}
	&& f_{n}\otimes _{1} \varphi_{n} \\
	 &=& \frac{ 4}{(n-1) ^{2}}\sum_{1\leq i<j\leq n} \sum_{1\leq k<l\leq n} \left( (h_{i} \widetilde{\otimes} h_{j})\otimes _{1} (g_{k} \widetilde{\otimes } g_{l} )\right) \nonumber\\
	&=& \frac{ 1}{(n-1) ^{2}}\sum_{1\leq i<j\leq n} \sum_{1\leq k<l\leq n}\int_{\mathbb{R}_+} \left[ h_{i}(x) h_{j}(u) + h_{j}(x) h_{i}(u)\right] \left[ g_{k}(y) g_{l}(u)+ g_{l}(y)g_{k}(u)\right] \textrm{d} u \nonumber\\
	&=& \frac{ 1}{(n-1) ^{2}}\sum_{1\leq i<j\leq n} \sum_{1\leq k<l\leq n}\left(  h_{i} \otimes g_{k} \mathbf{1}_{\{ j=l \}}+ h_{i} \otimes g_{l} \mathbf{1}_{\{ j=k \}}+ h_{j} \otimes g_{k} \mathbf{1}_{ \{ i=l \}} + h_{j}\otimes g_{l} \mathbf{1}_{ \{ i=k \}}\right)\nonumber \\
	&=&  \frac{ 1}{(n-1) ^{2}}\left[ \sum_{\substack{1 \leqslant i<j \leqslant n \\ 1 \leqslant k < j \leqslant n}} h_{i}\otimes g_{k} + \sum_{1\leq i<j<l\leq n }h_{i}\otimes g_{l} + \sum_{1\leq k<i<j\leq n} h_{j}\otimes g_{k} + \sum _{\substack{ 1\leq i<j\leq n \\ 1\leq i<l\leq n }} h_{j}\otimes g_{l} \right]\label{30i-8}
\end{eqnarray}
and for $n$ large enough,
\begin{eqnarray*}
	\left| 	f_{n}\otimes _{1} \varphi_{n} \right| &\leq & \frac{C}{n} \left| \sum _{1\leq i,k\leq n}h_{i} \otimes g_{k} \right|
\end{eqnarray*}
and then
\begin{eqnarray*}
	\norm{ f_{n}\otimes _{1} \varphi_{n}}{H ^{\otimes 2}}^{2} & \leqslant & \frac{C}{n^2} \sum_{1\leq i,k, i', k'\leq n} \left\langle h_{i}\otimes g_{k}, h_{ i'}\otimes g_{k '} \right\rangle_{ H ^{\otimes 2}}\\
	&=& \frac{C}{n^2} \sum_{1\leq i,k\leq n} \frac{1}{i^a k^a} \leq \frac{C}{n^{2a}} \xrightarrow[n \to \infty]{} 0.
\end{eqnarray*}
Consequently, from (\ref{30i-7}) and (\ref{30i-8}) it follows that (\ref{cdn}) is verified and $\left( (U_{n}, V_{n}), n\geq 1\right) $ converges in law, as $n\to \infty$, to $ F(1)\otimes F(1)$ and
\begin{equation}
	\label{30i-9}
	\E \left[ \left\langle \textrm{D}(-L) ^{-1} U_{n}, \textrm{D} V_{n} \right\rangle^{2} \right] \leqslant \frac{C}{n ^{2a}}.
\end{equation}
By Theorem \ref{thmchaos}, 
\begin{eqnarray*}
	d_{2} \left( \P_{(U_{n}, V_{n})}, F(1)\otimes F(1)\right) \leqslant &  & C  \E \left[ \left( 2(U_{n}+1) - \left\langle \textrm{D}(-L) ^{-1} U_{n}, \textrm{D}U_{n} \right\rangle _{H} \right)^2 \right]^{\frac{1}{2}}\\
	& + & C \E \left[ \langle \textrm{D} (-L) ^{-1}U_{n}, \textrm{D} U_{n}\rangle _{H} ^{2} \right]^{\frac{1}{2}} + d_{2} (\P_{V_{n}}, F(1)).
\end{eqnarray*}
The estimates (\ref{30i-1}), (\ref{30i-2}) and (\ref{30i-9}) give, for $n$ large,
\begin{equation*}
		d_{2} \left( \P_{(U_{n}, V_{n})}, F(1)\otimes F(1)\right) \leq C\left( n ^{-\frac{1}{2}}+ n ^{-a}\right). 
\end{equation*}

\section{Appendix}
Here we present the basics of Malliavin calculus and the proof of Lemma \ref{ll3}. 

\subsection{ Wiener-Chaos and Malliavin derivatives}

Here we describe the elements from stochastic analysis that we will need in the paper. Consider $H$ a real separable Hilbert space and $(W(h), h \in H)$ an isonormal Gaussian process on a probability space $(\Omega, {\cal{A}}, \P)$, which is a centered Gaussian family of random variables such that ${\bf E}\left[ W(\varphi) W(\psi) \right]  = \langle\varphi, \psi\rangle_{H}$. Denote by  $I_{n}$ the multiple stochastic integral with respect to
$W$ (see \cite{N}). This mapping $I_{n}$ is actually an isometry between the Hilbert space $H^{\odot n}$(symmetric tensor product) equipped with the scaled norm $\frac{1}{\sqrt{n!}}\Vert\cdot\Vert_{H^{\otimes n}}$ and the Wiener chaos of order $n$ which is defined as the closed linear span of the random variables $H_{n}(W(h))$ where $h \in H, \|h\|_{H}=1$ and $H_{n}$ is the Hermite polynomial of degree $n \in {\mathbb N}$
\begin{equation*}
	\forall x \in \mathbb{R}, \hskip0.5cm  H_{n}(x)=\frac{(-1)^{n}}{n!} \exp \left( \frac{x^{2}}{2} \right)
	\frac{\text{d}^{n}}{\text{d} x^{n}}\left( \exp \left( -\frac{x^{2}}{2}\right)
	\right).
\end{equation*}
The isometry of multiple integrals can be written as follows: for $m,n$ positive integers,
\begin{eqnarray}
	\mathbf{E}\left[ I_{n}(f) I_{m}(g) \right] &=& n! \langle \tilde{f},\tilde{g}\rangle _{H^{\otimes n}} \quad \mbox{if } m=n,\nonumber \\
	\mathbf{E} \left[ I_{n}(f) I_{m}(g) \right] &= & 0 \quad \mbox{if } m\not=n.\label{iso}
\end{eqnarray}
It also holds that
\begin{equation*}
	I_{n}(f) = I_{n}\big( \tilde{f}\big)
\end{equation*}
where $\tilde{f} $ denotes the symmetrization of $f$ defined by the formula 

$$\tilde{f}
(x_{1}, \ldots , x_{n}) =\frac{1}{n!} \sum_{\sigma \in {\cal S}_{n}}
f(x_{\sigma (1) }, \ldots , x_{\sigma (n) } ).$$

We recall that any square integrable random variable which is measurable with respect to the $\sigma$-algebra generated by $W$ can be expanded into an orthogonal sum of multiple stochastic integrals
\begin{equation}
	\label{sum1} F=\sum_{n=0}^\infty I_{n}(f_{n})
\end{equation}
where $f_{n}\in H^{\odot n}$ are (uniquely determined)
symmetric functions and $I_{0}(f_{0})=\mathbf{E}\left[  F\right]$.

\noindent Let $L$ be the Ornstein-Uhlenbeck operator
\begin{equation*}
	LF=-\sum_{n\geq 0} nI_{n}(f_{n})
\end{equation*}
if $F$ is given by (\ref{sum1}) and it is such that $\sum_{n=1}^{\infty} n^{2}n! \Vert f_{n} \Vert ^{2} _{{\cal{H}}^{\otimes n}}<\infty$. The pseudo-inverse of $L$, denoted $(-L)^{-1}$, is given by,  
\begin{equation*}
	(-L) ^{-1}(F)= \sum_{n\geq 1} \frac{1}{n}I_{n}(f_{n})
\end{equation*}
for $F$ as in (\ref{sum1}) with $\E[F]=0$.

\noindent For $p>1$ and $\alpha \in \mathbb{R}$ we introduce the Sobolev-Watanabe space $\mathbb{D}^{\alpha ,p }$  as the closure of
the set of polynomial random variables with respect to the norm
\begin{equation*}
	\Vert F\Vert _{\alpha , p} =\Vert (I -L) ^{\frac{\alpha }{2}} F \Vert_{L^{p} (\Omega )}
\end{equation*}
where $I$ represents the identity. We denote by $\text{D}$ the Malliavin  derivative operator that acts on smooth functions of the form $F=g(W(h_1), \dots , W(h_n))$ ($g$ is a smooth function with compact support and $h_i \in H$)
\begin{equation*}
	\text{D}F=\sum_{i=1}^{n}\frac{\partial g}{\partial x_{i}}(W(h_1), \ldots , W(h_n)) h_{i}.
\end{equation*}
The operator $D$ is continuous from $\mathbb{D}^{\alpha , p} $ into $\mathbb{D} ^{\alpha -1, p} \left( H\right).$

We will intensively use the product formula for multiple integrals.
It is well-known that for $f\in H^{\odot n}$ and $g\in H^{\odot m}$
\begin{equation}\label{prod}
	I_n(f)I_m(g)= \sum _{r=0}^{n\wedge m} r! \binom{n}{r} \binom{m}{r} I_{m+n-2r}(f\otimes _r g)
\end{equation}
where $f\otimes _r g$ means the $r$-contraction of $f$ and $g$ (see e.g. Section 1.1.2 in \cite{N}).

\noindent We also need to introduce the Skorohod integral integral (or the divergence operator), denoted by $\delta$,  which is the adjoint operator of $\textrm{D}$. Its domain is 
\begin{equation*}
	\text{Dom} (\delta) = \left\{ u \in L ^{2} \left(\Omega; H\right), \mathbf{E} \left[ \left| \langle \textrm{D}F, u \rangle_{H} \right| \right] \leq C \norm{F}{L^2 (\Omega)} \right\}  
\end{equation*}
and we have the duality relationship
\begin{equation}\label{dua}
	 \forall F \in \mathcal{S}, \forall u \in \text{Dom} (\delta), \hskip0.3cm \mathbf{E} \left[ F \delta (u) \right] = \mathbf{E} \left[ \langle \text{D}F, u \rangle _{H} \right]. 
\end{equation}

\subsection{Proof of Lemma \ref{ll3}}

\noindent {\bf Proof: } Let us first prove point 1. Concerning  the denominators in (\ref{fh}), notice that their new expressions come from the following equality easily derived by differentiating the left sides (and by being  careful with the bounds which change signs) :

\begin{equation}\label{23i-3}
	\left\{ \begin{array}{rcll}
		2(x+ \nu) p_{\nu} (x) & = &  - \int_{-\nu}^{x} u p_{\nu} (u) \ \textrm{d}u & \text{ on } (-\nu, + \infty)  \\
		2(x+ \nu) q_{\nu} (x) & = & \int_{x}^{- \nu} u q_{\nu} (u) \ \textrm{d}u  & \text{ on } (- \infty, - \nu). 
	\end{array} \right. 
\end{equation}
We will only focus on the numerator part. 

\noindent $\bullet$ \textit{Assume $x > -\nu$.} Remember that we have, for every ${\bf y}\in \mathbb{R} ^{n}$,

$$ f_h (x, \mathbf{y}) = \frac{1}{(x+\nu) p_{\nu} (x)} \int_{- \nu}^x p_{\nu} (u) \tilde{h}(u, \mathbf{y}) \ \textrm{d}u, $$
where $\tilde{h}(u, \mathbf{y}) = h(u,\mathbf{y}) - \E[h(Z_{\nu}, \mathbf{y})]$.  We will prove the expression (\ref{fhfub+}) in two steps. First, we will use Fubini's theorem to have an expression of $\tilde{h}$ in terms of $\frac{\partial h}{\partial x}$ and then we use  again  Fubini's theorem to get the desired formula for  $f_h$. 

Let us Fubini's theorem for $\tilde{h}$ by writing the expectation as :

$$ \E[h(Z_{\nu}, \mathbf{y})] = \int_{- \nu}^{+ \infty} h(z, \mathbf{y}) p_{\nu} (z) \ \textrm{d}z. $$

\noindent Hence, we have for every $u > - \nu$, we can write
\begin{eqnarray*}
	\tilde{h}(u, \mathbf{y}) = && \int_{-\nu}^{+ \infty} p_{\nu} (z) (h(u, \mathbf{y}) - h(z, \mathbf{y})) \ \textrm{d}z = \int_{-\nu}^{+ \infty} p_{\nu} (z) \left[ \int_z^u \frac{\partial h}{\partial x} (w, \mathbf{y}) \ \textrm{d}w \right] \ \textrm{d}z \\ 
	= & &  \int_{-\nu}^{u} p_{\nu} (z) \left[ \int_z^u \frac{\partial h}{\partial x} (w, \mathbf{y}) \ \textrm{d}w \right] \ \textrm{d}z - \int_{u}^{+ \infty} p_{\nu} (z) \left[ \int_u^z \frac{\partial h}{\partial x} (w, \mathbf{y}) \ \textrm{d}w \right] \ \textrm{d}z. 
\end{eqnarray*}
The last line allows us to switch the integrals. We get :

\begin{eqnarray}
	\tilde{h}(u, \mathbf{y}) = & &  \int_{- \nu}^{u} \frac{\partial h}{\partial x} (w, \mathbf{y})  \left[ \int_{- \nu}^w p_{\nu} (z) \ \textrm{d}z \right]  \ \textrm{d}w -  \int_{u}^{+ \infty} \frac{\partial h}{\partial x} (w, \mathbf{y})  \left[ \int_{w}^{+ \infty} p_{\nu} (z) \ \textrm{d}z \right]  \ \textrm{d}w\nonumber \\ 
	= & &  \int_{- \nu}^{u} \frac{\partial h}{\partial x} (w, \mathbf{y}) F_{\nu} (w) \ \textrm{d}w - \int_{u}^{+ \infty} \frac{\partial h}{\partial x} (w, \mathbf{y})  (1 - F_{\nu} (w))  \ \textrm{d}w.\label{new1}  
\end{eqnarray}
We now plug this last equality into $(\ref{fh})$ to have our expression.

\begin{eqnarray*}
	(x+\nu) p_{\nu} (x) f_h (x, \mathbf{y}) = &  & \int_{- \nu}^x p_{\nu} (u)   \int_{- \nu}^{u} \frac{\partial h}{\partial x} (w, \mathbf{y}) F_{\nu} (w) \ \textrm{d}w   \ \textrm{d}u \\ 
	& - & \int_{- \nu}^x p_{\nu} (u)  \int_{u}^{+ \infty} \frac{\partial h}{\partial x} (w, \mathbf{y})  (1 - F_{\nu} (w))  \ \textrm{d}w   \ \textrm{d}u.
\end{eqnarray*}
The first term can have its integrals switched, whereas for the second one, we need to write 

$$ \int_{- \nu}^x \int_{u}^{+ \infty} = \int_{- \nu}^x \int_{u}^x + \int_{- \nu}^x \int_{x}^{+ \infty}. $$
We obtain 
\begin{eqnarray*}
	(x+\nu) p_{\nu} (x) f_h (x, \mathbf{y}) = &  & \int_{- \nu}^x   \frac{\partial h}{\partial x} (w, \mathbf{y}) F_{\nu} (w) \int_{w}^{x} p_{\nu} (u) \ \textrm{d}u  \ \textrm{d}w    \\ 
	& - & \int_{- \nu}^x \frac{\partial h}{\partial x} (w, \mathbf{y})  (1 - F_{\nu} (w))   \int_{- \nu}^{w} p_{\nu} (u)  \ \textrm{d}u \ \textrm{d}w \\ 
	& - & \int_{x}^{+ \infty} \frac{\partial h}{\partial x} (w, \mathbf{y})  (1 - F_{\nu} (w))   \int_{- \nu}^{x} p_{\nu} (u)  \ \textrm{d}u \ \textrm{d}w.
\end{eqnarray*}
In terms of $F_{\nu}$, it writes 

\begin{eqnarray}
	&&	(x+\nu) p_{\nu} (x) f_h (x, \mathbf{y})\nonumber \\
	= &  & \int_{- \nu}^x   \frac{\partial h}{\partial x} (w, \mathbf{y}) F_{\nu} (w) (F_{\nu}(x) - F_{\nu} (w))  \ \textrm{d}w    - \int_{- \nu}^x \frac{\partial h}{\partial x} (w, \mathbf{y})  (1 - F_{\nu} (w))   F_{\nu} (w) \ \textrm{d}w \nonumber \\ 
	& - & \int_{x}^{+ \infty} \frac{\partial h}{\partial x} (w, \mathbf{y})  (1 - F_{\nu} (w))  F_{\nu} (x) \ \textrm{d}w \nonumber \\
	= & - & (1 - F_{\nu} (x)) \int_{- \nu}^x \frac{\partial h}{\partial x} (w, \mathbf{y})  F_{\nu} (w) \ \textrm{d}w -  F_{\nu} (x) \int_x^{+ \infty} \frac{\partial h}{\partial x} (w, \mathbf{y}) (1 - F_{\nu} (w)) \ \textrm{d}w. \label{23i-4}
\end{eqnarray}
By combining (\ref{23i-3}) and (\ref{23i-4}), we obtain the expression (\ref{fhfub+}).

\vskip0.2cm
\noindent $\bullet$ \textit{Assume $x < -\nu$.}  From (\ref{fh}), we have 

$$ f_h (x, \mathbf{y}) = \frac{-1}{(x+\nu) q_{\nu} (x)} \int_{x}^{- \nu} q_{\nu} (u) \tilde{h}(u) \ \textrm{d}u. $$
Since the integral is on $(x, - \nu)$, we cannot use the expression of $\tilde{h}$ derived on $(- \nu, + \infty)$. In fact, we need to compute it again, on $(- \infty, - \nu)$ this time. 

Let $u < - \nu$. We write
\begin{eqnarray*}
	\tilde{h}(u, \mathbf{y}) = & - &  \int_{-\nu}^{+ \infty} p_{\nu} (z) \left[ \int_u^z \frac{\partial h}{\partial x} (w, \mathbf{y}) \ \textrm{d}w \right] \ \textrm{d}z. \\
\end{eqnarray*}

\noindent By writing $\int_{- \nu}^{+ \infty} \int_{u}^z = \int_{- \nu}^{+ \infty} \int_{u}^{- \nu} + \int_{- \nu}^{+ \infty} \int_{- \nu}^z$, we use Fubini's theorem to have 

\begin{eqnarray*}
	\tilde{h}(u, \mathbf{y}) = & - &  \int_{u}^{- \nu} \frac{\partial h}{\partial x} (w, \mathbf{y})  \left[ \int_{- \nu}^{+ \infty} p_{\nu} (z)  \ \textrm{d}z \right] \ \textrm{d}w - \int_{- \nu}^{+ \infty} \frac{\partial h}{\partial x} (w, \mathbf{y})  \left[ \int_{w}^{+ \infty} p_{\nu} (z)  \ \textrm{d}z \right] \ \textrm{d}w \\
	= & - &  \int_{u}^{- \nu} \frac{\partial h}{\partial x} (w, \mathbf{y}) \ \textrm{d}w - \int_{- \nu}^{+ \infty} \frac{\partial h}{\partial x} (w, \mathbf{y}) (1 - F_{\nu} (w)) \ \textrm{d}w.
\end{eqnarray*}
And again, we plug this expression into the expression (\ref{fh}) of $f_h$, we find 

\begin{eqnarray*}
	- (x+\nu) q_{\nu} (x) f_h (x, \mathbf{y}) = & - &  \int_{x}^{- \nu} q_{\nu} (u)  \int_{u}^{- \nu} \frac{\partial h}{\partial x} (w, \mathbf{y}) \ \textrm{d}w \ \textrm{d}u \\
	& - & \int_{x}^{- \nu} q_{\nu} (u) \int_{- \nu}^{+ \infty} \frac{\partial h}{\partial x} (w, \mathbf{y}) (1 - F_{\nu} (w)) \ \textrm{d}w \ \textrm{d}u.
\end{eqnarray*}

\noindent The minus signs cancel. For both terms, both integrals can be switched by Fubini's theorem without using any Chasles relation this time. We get 

\begin{eqnarray*}
	&&	(x+\nu) q_{\nu} (x) f_h (x, \mathbf{y})\\
	= &  &  \int_{x}^{- \nu}  \frac{\partial h}{\partial x} (w, \mathbf{y}) \int_{x}^{w} q_{\nu} (u) \ \textrm{d}u \ \textrm{d}w + \left( \int_{x}^{- \nu} q_{\nu} (u) \ \textrm{d}u \right) \left( \int_{- \nu}^{+ \infty} \frac{\partial h}{\partial x} (w, \mathbf{y}) (1 - F_{\nu} (w)) \ \textrm{d}w \right) \\
	= & &   \int_{x}^{- \nu}  \frac{\partial h}{\partial x} (w, \mathbf{y}) \left[ \tilde{F}_{\nu}(x) - \tilde{F}_{\nu} (w) \right] \ \textrm{d}w +  \tilde{F}_{\nu} (x) \int_{- \nu}^{+ \infty} \frac{\partial h}{\partial x} (w, \mathbf{y}) (1 - F_{\nu} (w)) \ \textrm{d}w,
\end{eqnarray*}
so \ref{fhfub-}) is proven. 

  We now prove  2.  (actually  the proof is slightly adapted from \cite{D}).  Let $x \in \R$. Then

\begin{eqnarray*}
	\int_{- \infty}^x F_{\nu} (s) \ \textrm{d}s = & & \int_{- \nu}^x \int_{- \nu}^s p_{\nu} (u) \ \textrm{d}u \ \textrm{d}s = \int_{- \nu}^x p_{\nu} (u) \int_u^x \textrm{d}s \ \textrm{d}u \\ 
	= & & x F_{\nu} (x) - \int_{- \nu}^x u p_{\nu} (u) \ \textrm{d}u.
\end{eqnarray*}
\noindent We proceed in a similar way  for the second equality and we get 

\begin{eqnarray*}
	\int_x^{+ \infty} (1 - F_{\nu} (s)) \ \textrm{d}s = & & \int_x^{+ \infty} \int_{s}^{+ \infty} p_{\nu} (u) \ \textrm{d}u \ \textrm{d}s = \int_x^{+ \infty} p_{\nu} (u) \int_x^u \textrm{d}w \ \textrm{d}u \\ 
	= & & \int_x^{+ \infty} u p_{\nu} (u) - x(1 - F_{\nu} (x)).
\end{eqnarray*}
Concerning the point 3. in the statement,  we use the same idea as before,

\begin{eqnarray*}
	\int_x^{- \nu} \tilde{F}_{\nu} (w) \ \textrm{d}w = & & \int_x^{- \nu} \int_w^{- \nu} q_{\nu} (u) \ \textrm{d}u \ \textrm{d}w =\int_x^{- \nu} q_{\nu} (u) \int_x^u \textrm{d}w \ \textrm{d}u \\ 
	= & & \int_x^{- \nu} u q_{\nu} (u) \ \textrm{d}u - x \tilde{F}_{\nu} (x).
\end{eqnarray*}\qed

\end{document}